
\input eplain
\input amstex
\documentstyle{amsppt}
\loadbold


\chardef\oldatcatcode=\the\catcode`\@
\catcode `\@=11
\def\logo@{}
\catcode `@=\oldatcatcode


\NoBlackBoxes

\def\Z{\Bbb Z}

\def\C{\Bbb C}
\def\Q{\Bbb Q}
\def\h{\Bbb H}
\def\Qbar{\overline{\Bbb Q}}
\def\R{\Bbb R}

\def\p{\Bbb P}

\def\Vol{\operatorname{vol}}
\def\Gal{\operatorname{Gal}}
\def\Aut{\operatorname{Aut}}
\def\Isom{\operatorname{Isom}}
\def\id{\operatorname{id}}
\def\Li{\operatorname{Li}}
\def\D{\operatorname{D}}
\def\m{\operatorname{m}}
\def\Disc{\operatorname{\Delta}}
\def\disc{\operatorname{disc}}
\def\Res{\operatorname{Res}}
\def\PSL{\operatorname{PSL}}
\def\K{\operatorname{K}}

\def\A{\Cal A}
\def\c{\Cal C}

\def\M{\Cal M}

\def\o{\Cal O}

\def\B{ B}

\def\r{\overline r}

\def\a{\alpha}
\def\b{\beta}

\def\z{\zeta}

\def\IM{\text{Im}\,}


\def\maps{\colon\thinspace}
\def\Isom{\operatorname{Isom}}
\def\tr{\operatorname{tr}}
\def\PSLn#1#2{\PSL_{#1} #2}
\def\SL#1#2{\operatorname{SL}_{#1} #2}
\def\SnapPea{SnapPea}
\def\Snap{Snap}
\def\spandef#1#2{{  \left\langle  {#1}  \ \left| \   {#2} \right. \right\rangle }}


\def\Aa{9}
\def\Ab{10}
\def\Ac{11}

\def\ShapeEqn{35}
\def\ApolyOfN{36}


\def\mbox{\text}
\def\emph#1{{\it #1\/}}
\def\mathscr{\Cal}

\topmatter

\title{
Mahler's Measure and the Dilogarithm (II)
}\endtitle

\author David W. Boyd and Fernando Rodriguez-Villegas \
        with an appendix by Nathan M. Dunfield
\endauthor

\leftheadtext\nofrills{D. W. Boyd and F. Rodriguez-Villegas with N. M. Dunfield}

\address
Department of Mathematics,
University of British Columbia,
Vancouver, B.C., Canada  V6T 1Z2
\endaddress
\email  boyd\@math.ubc.ca  \endemail

\address 
 Department of Mathematics,
University of Texas at Austin,
Austin, TX 78712, USA
\endaddress
\email villegas\@math.utexas.edu \endemail

\address
Mathematics 253-37, 
Caltech, 
Pasadena, CA 91125, USA
\endaddress
\email dunfield\@caltech.edu \endemail


\date{
September 2005
}\enddate

\thanks The first author was sponsored in part by a grant from
NSERC. The second author was supported in part by a grant from the
NSF; he would also like to thank the Department of Mathematics at
Harvard University, where part of this work was done, for its
hospitality.
\endthanks

\abstract We continue to investigate the relation between the Mahler
 measure of certain two variable polynomials, the values of the
 Bloch--Wigner dilogarithm $D(z)$ and the values $\zeta_F(2)$ of zeta
 functions of number fields. Specifically, we define a class $\A$ of
 polynomials $A$ with the property that $\pi\m(A)$ is a linear
 combination of values $D$ at algebraic arguments. For many
 polynomials in this class the corresponding argument of $D$ is in the
 Bloch group, which leads to formulas expressing $\pi \m(A)$ as a linear
 combination with unspecified rational coefficients of $V_F$ for
 certain number fields $F$ ($V_F := c_F\zeta_F(2)$ with $c_F>0$ an
 explicit simple constant). The class $\A$ contains the
 $A$-polynomials of cusped hyperbolic manifolds. The connection
 with hyperbolic geometry often provides means to prove identities of
 the form $\pi \m(A)= r V_F$ with an explicit value of $r\in \Q^*$. We
 give one such example in detail in the body of the paper and in the
 appendix. 
\endabstract

\endtopmatter

\document


\baselineskip=14pt

\vskip.5cm
\subheading {Introduction}
The main goal of this paper is to understand and expand on Smyth's
remarkable discovery
$$
\m(x+y+1)=L'(\chi,-1).
$$

This formula relates the Mahler measure of a polynomial to the value
at  $s = 2$ of the zeta function of a number field.  Essentially all
known examples of this kind of formula arise from the combination of
the following three facts.

\vskip.5cm {\bf A.} There exists a class $\A$ of two variable
Laurent polynomials $A$ with the property that $\pi \m(A)$ equals an
explicit rational linear combination of values of the Bloch--Wigner
dilogarithm at algebraic arguments.

\vskip.5cm {\bf B.} The volume of hyperbolic 3-manifolds  can be
expressed as a rational linear combination of the Bloch--Wigner
dilogarithm at algebraic arguments.

\vskip.5cm {\bf C.} The volume of an arithmetic hyperbolic 3-manifold
$N$ can be expressed in terms of $\zeta_F(2)$, where $\zeta_F$ is the
zeta function of an associated number field $F$.

\vskip.5cm Let us make some remarks about these facts: {\bf A} is
analogous to the fact that for all one-variable Laurent polynomials
$A$, $\m(A)$ can be expressed in terms of $\log|\alpha|$, where
$\alpha$ runs through the roots of $A$.  Roughly speaking, the class
$\A$ consists of the minimal polynomials of pairs of rational
functions $x,y$ on some algebraic curve $X$ whose symbol $\{x,y\}$ is
trivial in $K_2(X)$.  {\bf B} is a consequence of the fact that a
hyperbolic 3-manifold can be decomposed into hyperbolic tetrahedra, 
whose volumes can be expressed in terms of the Bloch--Wigner
dilogarithm.  {\bf C} is due to Borel in general and to Humbert for
imaginary quadratic fields.

These facts are also interrelated. There is an invariant $A(x,y)\in
\Q[x,x^{-1},y,y^{-1}]$ of a 1-cusped hyperbolic manifold $M$ (for
example the complement of a hyperbolic knot) called the {\it
$A$-polynomial} of $M$ (or of the knot). It is proved in \cite{CCGLS}
that $A\in \A(\Q)$. In some cases, but not always, this leads to identities of
the form
$$
\pi \m(A)=\Vol(M).
\tag{a}
$$

A theorem of Borel and others (see \S 4) guarantees that for number
fields $F$ with only one complex embedding certain rational linear
combinations of the Bloch--Wigner dilogarithm at algebraic arguments
is commensurable with $V_F$, an invariant of $F$ which equals a simple
factor times $\zeta_F(2)$ \cite{Bor,Za}. (We say that two $a,b\in\C$
are {\it commensurable} and write $a\sim b$ if there exists a non-zero
$r\in \Q$ such that $a=rb$.) Combining this theorem with {\bf A}, we
find for many polynomials $A\in \A(\Q)$ that
$$
\pi \m(A) \sim V_F
\tag{b}
$$
for some number field $F$.

Unfortunately, Borel's theorem does not provide an a priori bound on
the height of the rational number implicit in (b) and hence there is
no way to way to pin down this number by a numerical
calculation. However, the most plausible numerical value obtained
by numerical calculations to high accuracy is usually a rational
of very small height, suggesting that these predicted values are likely
correct.

Apart from very simple examples, the only way we see how to obtain an
identity
$$
\pi \m(A)= r V_F
\tag{c}
$$ with a provable value of $r\in \Q^*$ is to combine all three facts:
First, obtain an identity of type (a) by means of {\bf A} and {\bf B};
then, relate the manifold $M$ explicitly with an arithmetic manifold
$N$ and use {\bf C}. We consider in detail  one such example for which
we give two different proofs (see \S 6 example 3 and the appendix).

 From our point of view, it would be very useful to have an a priori
 upper bound for the height of the rational number in Borel's
 theorem. Equivalently, it would be good to know that an analogue of
 Lehmer's conjecture for $K_3$ of number fields is valid \cite{RV2}.

\subheading{ 1. Definition of the class $\A$}
\noindent

For a Laurent polynomial
$$
A=\sum_{(m,n) \in \Z^2}a_{(m,n)}\; x^m y^n \; \in \C[x^{\pm1},y^{\pm1}]
$$
we define its {\it Newton polygon} $\Delta \subset \R^2$ as the
convex hull of the finitely many points $(a,b) \in \Z^2$ for which
$a_{m,n}$ is non-zero. We call $a_v$, where $v$ is a vertex of
$\Delta$ a {\it vertex coefficient} of $A$.

We will usually give polynomials by their matrix of coefficients, whose
$(i,j)$ entry is $a_{(j,i)}$ (so that the coefficients of increasing
powers of $y$ are read top to bottom and those of $x$ left to right),
from which one easily determines $\Delta$. (We will also drop $0$
entries unless this could lead to confusion.)

Let $X$ be a smooth projective algebraic curve defined over $\C$ and
let $\C(X)$ be its function field. Let $x,y\in\C(X)$ be two
non-constant rational functions and let $S\subset X$ be the set of
zeros and poles of $x$ or $y$. The image of the rational map 
$(x,y): X\setminus S \longrightarrow  \C^\star \times \C^\star$
is of dimension $1$; let $A\in \C[x,y]$
be a defining equation. In terms of this equation the condition
$\{x,y\} \in K_2(X)\otimes \Q$ is equivalent to $A$ being {\it
tempered}, i.e., the roots of all the face polynomials $A$ are roots
of unity (see \cite{RV1}).  If this is the case, it is not hard to see
that we may scale $A$ so that it has all its vertex coefficients $a_v$
equal to roots of unity and, with some notational abuse, we will simply
call an $A$ so normalized the {\it minimal polynomial} of the pair
$x,y$.  In particular, if $X,x,y$ are all defined over $\Q$ we will
normalize $A$ so that its vertex coefficients are $\pm 1$, which
determines it uniquely up to sign.

In this paper we will be mostly concerned with pairs $x,y$ as above
satisfying
$$
\{x,y\} =0  \qquad  \text{in  } \K_2(X) \otimes \Q;  \tag {\bf K}
$$
or, equivalently,   for some $r_j\in \Q$ and $z_j\in \C(X)^*$
we have
$$
x\wedge y= r_1\langle z_1\rangle + \cdots + r_n\langle z_n\rangle
\tag {\bf T}
$$
in $\wedge^2(\C(X)^*)\otimes\Q$, where to simplify the notation we
have set $\langle z\rangle:=z\wedge(1-z)$.  We will refer to such an
identity as a {\it triangulation} of $x\wedge y$. The reason for this
terminology is that this identity is an algebraic analogue of
decomposing a hyperbolic 3-manifold into ideal tetrahedra.

For a field $F\subset \C$ we let $\A(F)$ be the collection of Laurent
polynomials $A$ such that, up to multiplication by some monomial
$x^ry^s$, $A$ is the minimal polynomial of a pair $x,y$ satisfying
({\bf K}) with $X,x,y$ defined over $F$.

A polynomial $A(x,y)$ is called {\it reciprocal} if
$A(x^{-1},y^{-1})=\pm x^ay^bA(x,y)$ for some integers $a,b$.  For
example, the $A$-polynomial of a $1$-cusped manifold is
reciprocal. For reciprocal polynomials the map $\sigma: (x,y) \mapsto
(x^{-1},y^{-1})$ extends to an involution of the smooth completion $X$
of the zero locus of $A$. This involution fixes the symbol $\{x,y\}$
as
$$
\{x,y\}^\sigma=\{x^{-1},y^{-1}\}=\{x,y\}
$$
by the bi-multiplicativity of symbols in $K_2$.

 For example, if $X$ is an elliptic curve and $\sigma$ has fixed
 points then $X/\langle \sigma \rangle$ is isomorphic to $\p^1$. Hence
 if $\{x,y\} \in K_2(C)$ then by Galois descent $\{x,y\}$ vanishes
 modulo torsion (see example 8 below for an application of this
 remark).

 In general, the involution on a generic elliptic curve $X/\Qbar$
 (i.e. for which $\Aut(X)$ is just $\pm \id_X$) is of the form
$$
\sigma: \quad P\mapsto \epsilon P +Q,
$$ where $\epsilon=\pm 1$ and $Q\in X(\Qbar)$ satisfies $Q+\epsilon
Q=0$ (vacuous if $\epsilon =-1$). The involutions with $\epsilon=+1$
have no fixed points; they generate the center of
$\Isom(X)$. Involutions with $\epsilon=-1$ have four fixed points;
namely, those $P$ which satisfy $2P=Q$; their field of definition is
at most of degree $4$ over that of $Q$. Two involutions with
$\epsilon=-1$ are conjugate modulo $\Isom(X)$ if and only if the
corresponding points $Q$ are congruent modulo $2X(\Qbar)$.

We will say that a number field is of {\it type} $d, [r_1,r_2], D$ if is
of degree $d$ over $\Bbb Q$, has $r_1$ real and $r_2$ complex embeddings
into $\Bbb C$ and has discriminant $D$.

Many of our examples refer to the census of hyperbolic 3-manifolds
that can be triangulated by $7$ or fewer ideal tetrahedra \cite{HW,
CHW} and distributed with Week's program SnapPea \cite {W}.  We have
made extensive use of SnapPea and the related program Snap described
in \cite{CGHN} in our computations.  Significant use has been made of
the number theory system GP-PARI and the computer algebra systems
Maple and Macaulay2.  Although all numerical computations were
performed to at least 28 decimal place precision, in this paper we
exhibit at most 10 decimal places of the corresponding numbers.

\subheading{2. Mahler measure and the class $\A(F)$}
The {\it logarithmic Mahler measure} of a non-zero
Laurent polynomial  $A \in \C
[x_1^{\pm1}, \ldots , \;  x_n^{\pm1}]$
is defined  as
$$
\m(A) := \int_0^1 \cdots \int_0^1 \, \log \left| A(e^{2\pi i \theta_1}
\, , \ldots , \;
      e^{2\pi i\theta_n}) \right| d\theta_1 \cdots d\theta_n \;
\tag{1}
$$
and its {\it Mahler measure}  as $M(A)=e^{\m(A)}$,
 the geometric mean of $|A|$ on the torus
$$
T^n:=\{ (z_1,
\ldots , z_n) \in \C^n | \; |z_1|= \ldots = |z_n| = 1 \}.
$$

If $A(x) = a_0\prod_{j=1}^d(x-x_j) \in \C[x]$ is  a polynomial in one
variable then  Jensen's formula yields
$$
\m(A)= \log |a_0| + \sum_{j=1}^d \log^+ |x_j|,
\tag{2}
$$
where $\log^+|z|= \log |z|$ if $|z|\geq 1$ and  $0$ otherwise.

We will see in what follows under what conditions on $A$ a formula of
this nature, with $\log$ replaced by the dilogarithm, holds.

If $A(x,y) \in \C[x,y]$ is a polynomial in two variables we may think
of it as a polynomial in $x$ with coefficients which are polynomials
in $y$ and write
$$
A(x,y)=a_0(y)\prod_{j=1}^d(x-x_j(y))
\tag{3}
$$
where $x_j(y)$ are algebraic functions of $y$. Integrating the $x$
variable using Jensen's formula we obtain
$$
\m(A)= \m(a_0) + \sum_{j=1}^d \frac1{2\pi i} \int_{|y|=1} \log^+
|x_j(y)|\; \frac{dy}y,
\tag{4}
$$

We assume from now on that $A$ is absolutely irreducible. There is no
loss in generality since clearly $\m(A_1 A_2)=\m(A_1) \m(A_2)$. Let
$Y\subset \C^*\times \C^*$ be the zero locus of $A$ and let $X$ be a
smooth projective completion of $Y$.

 We want to rewrite the right hand side of (3) in a more manageable
way. To this end, regarding $x$ and $y$ as rational functions on $X$,
we define
$$
\eta(x,y):=\log|x|d\arg y -\log |y| d\arg x,
\tag{5}
$$
a  real differential $1$-form on $X\setminus S$, where $S$ is
the finite set of points of $X$ where either $x$ or $y$ have a zero or
a pole.

We have
$$
\m(A)=\m(a_0)+ \frac{1}{2\pi} \int_\gamma \eta(x,y),
\tag{6}
$$
where $\gamma$ is an oriented path on $X$ projecting to the
intersection of $Y$ with the set $|y|=1, |x|\geq 1$. It is not hard to
verify that we may chose $\gamma$ so that its boundary is
$$
\partial\, \gamma = \sum_k \epsilon_k\, [w_k], \qquad \epsilon_k = \pm 1
\tag{7}
$$ 
with $|x(w_k)|=|y(w_k)|=1$ for all $k$; moreover, the points $w_k$
are defined over $\Qbar$.

It is easy to check that
$$
d\eta =\IM \left(\frac{dx}x\wedge \frac{dy}y\right)
\tag{8}
$$
which vanishes since $\dim X=1$ and hence $\eta(x,y)$ is a closed
differential.

When does $\eta(x,y)$ extend to all of $X$? It is not hard to check that
for any point $w\in X(\C)$ we have
$$
\frac1{2\pi}\int_{C_w}\eta(x,y)=\log|(x,y)_w|,
\tag{9}
$$
where $C_w$ is any sufficiently small positively oriented circle with
center $w$ and
$$
(x,y)_w=(-1)^{w(x)w(y)} \left. \frac{x^{w(y)}}{y^{w(x)}}\right|_w
\tag{10}
$$
is the tame symbol at $w$. As a consequence, if $A$ is tempered then
$\eta$ extends to all of $X$ (see \cite{RV1}). 

The tame symbol gives rise to a homomorphism
$K_2(\C(X))\longrightarrow \C^*$ Hence, if $\{x,y\}$ satisfies ({\bf
K}) then all tame symbols $(x,y)_w$ are automatically torsion.

 A fundamental example is the following. Let $X=\p^1$, $x=t, y=1-t$,
where $t$ is a parameter on $X$, so that $A=x+y-1$. Since all tame
symbols are clearly torsion $\eta(t,1-t)$ extends to $\p^1$ and
is hence exact. Indeed, we have
$$
\eta(t,1-t) =d\D(t),
\tag{11}
$$
where $\D$ is the Bloch--Wigner dilogarithm
$$
\D(t)=\IM(\Li_2(t))+\arg(1-t)\log|t|
\tag{12}
$$

In general, if we have a triangulation
$
x\wedge y= r_1\langle z_1\rangle + \cdots + r_n\langle z_n\rangle
$
as in {\bf (T)}, then
$$
\eta(x,y)=dV,
\tag{13}
$$
where
$$
V=\D(\xi), \qquad \xi=r_1\left[z_1\right]+\cdots+r_n\left[z_n\right]
\tag{14}
$$
and $\D$  is extended to the group ring $\Z[\C(X)]$ by linearity.

We can now use Stokes theorem to conclude that
$$
2\pi\m(A)= \sum_k \epsilon_k V(w_k)
\tag{15}
$$
which we may also write as
$$
2\pi\m(A)= D(\xi)
\tag 16
$$
where
$$
\xi=\sum_k\epsilon_k \xi_k, \qquad \qquad \xi_k= \sum_j
\left[z_j(w_k)\right].
\tag 17
$$
In a sense, (15) is the analogue of (2) for two variable
polynomials in $\A$.

We summarize the above discussion in the following

\proclaim{Theorem 1}
 Let the  notation be as above and assume ({\bf K}) holds.   Then
$$
2\pi\m(A)= \sum_k \epsilon_k D(\xi_k)
$$
with $\epsilon_k=\pm 1, \xi_k \in \Z\left[\Qbar\right]$ as defined in (17).
\endproclaim

In particular, Theorem 1 applies to the A-polynomials of hyperbolic
3-manifolds.

Going back to the basic example $X=\p^1,x=t,y=1-t,A=x+y-1$, we take
$z_1=t$ and $\gamma$  the circle $1-e^{i\theta}$ with $ \theta \in
[-2\pi/3,2\pi/3]$ in the $t$-plane. We have
$$
\partial \gamma =\left[w_1\right]-\left[\overline w_1 \right],
$$
where $w_1 = (e^{\pi i/3},e^{2\pi i/3})$, and hence
$$
2\pi\m(x+y-1)=\D(e^{\pi i/3})-\D(e^{-\pi i/3}).
\tag 18
$$
Using the expansion
$$
\D(e^{i\theta})=\sum_{n=1}^\infty \frac{\sin(n\theta)}{n^2}, \qquad
\theta \in \R
$$
it is easy to check that the right hand side of (18) equals
$
(3\sqrt{3}/2)L(\chi,2),
$
with $\chi$ the quadratic character attached to $\Q(\sqrt{-3})$,
which yields Smyth's result
$$
\m(x+y+1)=L'(\chi,-1).
$$

\subheading{3. The five term relation for the dilogarithm}
To illustrate the identity (13) consider the affine variety $X$ defined by
following set of equations in the five variables $z_0,\ldots, z_4$.
$$
\aligned
1-z_1&=z_2z_0\\
1-z_2&=z_3z_1\\
1-z_3&=z_4z_2\\
1-z_4&=z_0z_3\\
1-z_0&=z_1z_4
\endaligned
\tag 19
$$
(they are obtained from the first by cyclically permuting the
variables with the indices read modulo $5$).  It is  a simple
matter to verify that
$$
z_0\wedge (1-z_0) + \cdots + z_4\wedge (1-z_4)=0
\tag 20
$$
in $\wedge^2(\C(X)^*)$. It follows that $V$ must be constant since
$dV$ vanishes; in fact this constant must be zero since, for example,
$\D(t)=0$ for $t \in \R$ and $(\phi,\ldots,\phi)\in X$,  where
$\phi=(-1+\sqrt{5})/2$. Therefore,
$$
\D\circ z_0+\cdots+\D\circ z_4=0, \qquad (z_0,\ldots,z_4) \in X.
\tag 21
$$
 This statement does actually have content since $X$ has positive
dimension; in fact, it is a surface birational to $\p^2$ (for example,
$z_0$ and $z_1$ determine the other variables uniquely giving a
birational isomorphism).  The projective completion $\bar X$ of $X$
was shown by Elkies to be a del Pezzo surface. The identity (21) is
the $5$-term relation satisfied by the dilogarithm.

\subheading{4. Mahler's measure and special values of Dedekind zeta functions} 
The relation between values of the
 dilogarithm and special values of $L$-functions we encountered at the
 end of section 2 is part of a more general phenomenon (due to Borel
 and others), which we now describe.

For any field $F$  and $\xi =\sum_j n_j \left[a_j\right] \in
\Z[F]$ we define
$$
\partial(\xi):=\sum_j n_j (a_j\wedge(1-a_j)) \in\wedge^2(F^*),
\tag 22
$$
where the corresponding term in the sum is omitted if $z_i=0,1$,
and let
$$
A(F)= \ker(\partial)
$$
Next we define the group
$$
C(F)= \left\{\left[x\right]+\left[y\right]+\left[\frac{1-x}{1-xy}\right]+
\left[1-xy\right]+\left[\frac{1-y}{1-xy}\right]
\;|\; x,y \in F,\; xy \neq 1\right\}\;.
$$
It is not hard to check that $C(F)\subset A(F)$. We define the {\it
Bloch group} as the quotient
$$
\B(F)=A(F)/C(F)\;.
$$
The $5$-term relation for the dilogarithm (21), guarantees that $D$
induces a well defined function on  $\B(\C)$ (still denoted by $D$).

Clearly an embedding of fields $\sigma: F \longrightarrow L$ extends by
linearity to a map $\sigma: \B(F) \longrightarrow \B(L)$ and in
particular  $\Gal(F/\Q)$ acts on $\B(F)$ if $F/\Q$ is Galois.

Given a number field $F$ denote by
$\Disc_F$ its  discriminant, $\zeta_F$ its zeta function and $r_1,r_2$
for the number of real and pairs of complex embeddings respectively.
The following results can be found in \cite{Za}:

\proclaim{Theorem A}
Let $F$ be a number field of degree $n$ with $r_2=1$ and let
 $\sigma$ be one its complex embeddings. Then $\B(F)$ is
of rank $1$ and if  $\xi =\sum_j n_j \left[a_j\right] \in
\Z[F\setminus \{0,1\}]$ represents a non-torsion element of $\B(F)$
then
$$
D(\sigma(\xi)) = r\;\frac{|\Disc_F|^{3/2}}{\pi^{2(n-1)}}\;\zeta_F(2)\;,
$$
for some $r \in \Q^*$.
\endproclaim
We will also need the following {\it Galois descent} property of the
Bloch group.
\proclaim{Theorem B} Let $L/F$ be a Galois extension of number
fields with $G=\Gal(L/F)$. Then
$$
\B(F)\otimes_\Z\Q=\B(L)^G\otimes_\Z\Q\;.
$$
\endproclaim
We want to apply these theorems in conjunction with (16) to obtain
relations between the Mahler measure of certain polynomials and
special values of $L$-functions analogous to Smyth's result. We find
that
$$
\partial(\xi_k)=\sum_j z_j(w_k)\wedge(1-z_j(w_k))=x(w_k)\wedge y(w_k)
\tag 23
$$

Hence we have the following

\proclaim { Theorem 2} Under the assumptions of Theorem 1, we have
$$
\xi_k \in\B(\Qbar) \qquad \text{if and only if} \qquad x(w_k)\wedge y(w_k)=0.
$$
\endproclaim

Given $\xi\in \Z[\Qbar]$ we will say that $\xi$ is {\it defined over
$F\subset \Qbar$} if every $\sigma \in\Gal(\Qbar/\Q)$ which satisfies
$\xi^\sigma=\xi$ is the identity on $F$.

Putting everything together we obtain the following.

\proclaim{Theorem 3}  
 Let the  notation be as above and assume ({\bf K}) holds so that
 the formula for $\pi m(A)$ of Theorem 1 holds.     Assume in addition
 that  he field of definition $F_k$ of $\xi_k$ has only one complex
embedding for all $k$  and that
$$
x(w_k)\wedge y(w_k)=0, \qquad \text{for all $k$}.
\tag {\bf B}
$$
so that $\xi_k \in\B(\Qbar)$ by Theorem 2.
Then
$$
\pi \m(A)= \sum_k r_k V_{F_k}
$$
for some rational numbers $r_k\in \Q$,
where
$$
V_F := {|\Disc_F|^{3/2}  \over (2\pi)^{2n-2}}\;\zeta_F(2)\;.
\tag 24
$$
\endproclaim

Note that in the paper \cite{BRV} we used instead the related quantity
$Z_F = 6V_F/\pi$.  Here $V_F$ seems more natural since it is the
volume of a hyperbolic orbifold $\Gamma_F\backslash\h^3$ \cite {Bor}.
For the quadratic field $F$ of discriminant $-f$, we have $V_F = \pi
d_f/6$, where $d_f = L'(\chi,-1)$, $\chi$ denoting the non-principal
character of conductor $f$.

\subheading{5. Curves of genus 0}

At this point it might seem unlikely that we will ever find examples
that satisfy all the conditions of Theorem 2 other than Smyth's example.
We claim that there are in fact many. We start with the simplest situation
of rational curves.

\example {Example 1 -- A rational curve}

Let $X=\p^1$ with $t$ a parameter as before and consider
$$
x=\frac{t^2+t+1}{(t-1)^2}, \qquad y= \frac{3t^2}{(t-1)^2}.
$$
It is not hard to check directly that all tame symbols $(x,y)_w$
are roots of unity. Alternatively, one can compute the minimal
polynomial relation between $x$ and $y$
$$
A=x^2-2xy-2x+1-y+y^2
$$
and check that it is tempered. Indeed, in matrix form $A$ is
$$
\matrix
\format \r\quad & \r\quad & \r\\\
1&-2&1\\
-1&-2&\\
1&&\\
\endmatrix
$$
with face polynomials $(u-1)^2,(u-1)^2,u^2-u+1$ all whose roots are
roots of unity.

On $\p^1$ any pair of non-zero rational functions $x,y$ with trivial
tame symbols $(x,y)_w$ at all points $w$ automatically satisfies {\bf (K)}
\cite {Mi}. Concretely, we may compute this decomposition (over $\Qbar$) by
means of the following identity (which we learned from J.~Tate)
$$
(t-a)\wedge(t-b)=\left(\frac{t-a}{b-a}\right)\wedge
\left(1-\frac{t-a}{b-a}\right) + (t-a)\wedge(a-b)-(t-b)\wedge(b-a)
\tag 25
$$
For our example this yields the following
$$
V(t)=\D(\xi(t))
$$
where
$$
\xi(t)=2\left[\frac{t-\zeta_3}{-\zeta_3}\right]
        +2\left[\frac{t-\zeta_3^{-1}}{-\zeta_3^{-1}}\right]
        -2\left[\frac{t-\zeta_3}{1-\zeta_3}\right]
        -2\left[\frac{t-\zeta_3^{-1}}{1-\zeta_3^{-1}}\right]
        -4[1-t]
$$
To find the points $w_k$ we can proceed as follows. Let
$$
A^*(x,y):=x^2y^2A(1/x,1/y)=(y^2-y+1)x^2-2(y^2+ y)x + y^2
$$
be the reciprocal of $A$. Any solution of $A(x,y)=0$ with
$|x(w_k)|=|y(w_k)|=1$ will also be a solution of $A^*(x,y)=0$. Hence,
the $x$-coordinates of the points on $A(x,y)=0$ with $|x|=|y|=1$ will
satisfy the equation
$$
\Res_y(A,A^*)=(x^2 - 4x + 1)^2(x^4 - 2x^3 -2x + 1)=0,
$$
where $\Res_y$ denotes the resultant in the variable $y$. Of the
roots of this equation only two (roots of the second factor) have
$|x|=1$, namely
$$
x_1=-0.3660254037 + 0.9306048591 i
$$
and its complex conjugate. We set $x = x_1$ in $A(x,y)=0$ and solve for
$y$ and find one value with $|y|=1$, namely
$$
y_1= -0.7320508075 - 0.6812500386 i.
$$
Solving for the corresponding value of $t$ we find
$$
t_1= 0.1339745962 + 0.5899798397 i
$$
These three numbers $x_1,y_1$ and $t_1$ lie in the splitting field
of $x^4 - 2x^3 - 2x + 1$, a number field $F$ of type $4, [2,1], -1728$.

It is not hard to verify that
$\partial(\gamma)=\left[w_1\right]-\left[\overline{w_1}\right]$, where
$w_1=(x_1,y_1)$; the correct sign for each point
equals the sign of the corresponding value of
$$
-\IM\left(\frac{y\frac{\partial P}{\partial y}}
{x\frac{\partial P}{\partial x}}\right).
$$

We conclude that
$$
2\pi \m(A)= \D(\xi(t_1))
$$
 The condition {\bf (B)} of Theorem 3 follows from the fact,
which is easy to verify, that $y_1=x_1^2$. That such a relation holds
is not as surprising as it might first appear if we notice that
both $x_1$ and $y_1$ lie in the same rank $1$ subgroup of $\o_F^*$,
the units of $F$, namely, the kernel of the homomorphism
$$
\aligned
\o_F^*&\qquad \longrightarrow \qquad  \C^*\\
u & \qquad \mapsto \qquad  \log|\sigma(u)|
\endaligned\;
$$
where $\sigma$ is either complex embedding of $F$.

Since $\xi(t_1)$  is stable
under the non-trivial Galois automorphism of $\Gal(F(\zeta_3)/F)$ it
is defined over $F$ and by Theorem 3.
$$
\pi \m(A)=r \; V_F
$$
for some $r\in \Q$ (which numerically equals $1/6$ to high precision).

\endexample

\example{Example 2 -- A second rational curve}

 A completely analogous example is the following.
Again we let $X=\p^1$ with $t$ a parameter and let
$$
x=\frac{2}{t^2-t+1}, \qquad y= \frac{t^2+t+1}{t^2-t+1}.
$$
Then
$$
A=x^2-xy+y^2-x-2y+1
$$
with matrix of coefficients
$$
\matrix
\format \r\quad & \r\quad & \r\\\
1&-1&1\\
-2&-1&\\
1&&\\
\endmatrix
$$
is tempered and
$$
V(t)=\D(\xi(t))
$$
with
$$
\xi(t)=-\left(\left[-t+\zeta_6\right]
        +\left[-t+\zeta_6^{-1}\right]
        +\left[\frac{\zeta_3}2(t-\zeta_6)\right]
        +\left[\frac{\zeta_3^{-1}}2(t-\zeta_6^{-1})\right]\right)
$$
The boundary of $\gamma$ again is of the form
$\partial(\gamma)=\left[w_1\right]-\left[\overline{w_1}\right]$ with
$w_1=(x_1,y_1)$ corresponding to $t=t_1=1.5174899135 i$ a root of
$t^4+t^2-3=0$. The number field $F$ is now the splitting field of
$x^4+x^2-3$, a field of type $4, [2,1], -507$.

Since again we have $y_1=x_1^2$ Theorem~3 applies and we obtain
$$
\pi \m(A)=r \; V_F
$$
for some $r\in \Q$ (which numerically equals $7/6$ to high precision).

\endexample

\subheading{6. Curves of positive genus}

 Given a curve $X$ of positive genus and two rational functions $x,y$
is not very easy to verify if ({\bf K}) holds and even if we know that
it does is not easy to find the actual triangulation ({\bf T}) of
$x\wedge y$.   One situation in which this is always possible is for
A-polynomials of hyperbolic manifolds and many, but not all, of our
examples are of this type.    It appears that the genus of the curve defined
by an A-polynomial may be arbitrarily large.   Once a triangulation is known
Theorem 1 applies to give a formula for $\pi m(A)$ as a sum of dilogarithms
of algebraic numbers.   However, Theorem 2 may not apply as Example 7
below illustrates.     
Our first example illustrates how a triangulation of
$x \wedge y$ may be found in case $A(x,y) = 0$ is an elliptic curve without knowing that
$A(x,y)$ is an A-polynomial

\example{Example 3 -- An elliptic curve} Our first example illustrates how a triangulation of
$x \wedge y$ may be found in case $A(x,y) = 0$ is an elliptic curve without knowing that
$A(x,y)$ is an A-polynomial.     The Appendix by Dunfield shows that in fact this
polynomial {\it is} an A-polynomial and leads to an independent evaluation of $m(A)$.

Let
$$
A=(x^2+x+1)(y^2+x)+3x(x+1)y
\tag 26
$$
with matrix of coefficients
$$
\matrix
\format \r\quad & \r\quad & \r\quad &\r\quad & \r\\
&1&1&1&\\
&3&3&&\\
1&1&1&&\\
\endmatrix
$$

The affine curve $Y\subset \C^*\times\C^*$ defined by the zero locus
of $A$ is singular, its only singularity is the ordinary double point
$(1,-1)$. A smooth projective completion of $Y$ is the elliptic curve
of conductor $15$ with minimal Weierstrass model
$$
X: \qquad v^2+uv+v=u^3+u^2.
$$
To see this we compute the discriminant of $A$ as a polynomial in $y$
and find
$$
\disc_y(A)=-(x-1)^2x(4x^2+7x+4)
$$
from which we deduce that $Y$ is birational to the elliptic curve
$$
w^2=-x(4x^2+7x+4)\;.
$$
A calculation then shows that
$$
\aligned
x&=-(u+1)\\
y&=\frac{v(v+1)}{u(u+v)}
\endaligned
$$
gives a birational map between $X$ and $Y$.

 The functions $x$ and $y$ on $X$ have divisors supported on the
 subgroup $H\subset X(\C)$ of order $8$ generated by
 $Q=(\zeta_3,-\zeta_3)$, where $\zeta_3\in\C$ is a fixed primitive
 third root of unity. We claim that there is a triangulation (13) for
 $x\wedge y$. In fact, we can find it using rational functions on $X$
 such that $f$ and $1-f$ have divisors supported on $H$.

Note that  $2Q=(0,0)$ is a rational point of order $4$. Let $\tau$
be the involution on the function field $\Q(X)$ determined by the map
$R\mapsto R+2Q$ for points $R$ of $X$. If $f,1-f$ is a pair functions
with divisors supported on $H$ so will $f^\tau, 1-f^\tau$.

 We find the pairs
$$
\aligned
-u,&\quad u+1\\
-v,&\quad v+1\\
-(u+v),&\quad u+v+1\\
 -\frac{v+1}u,&\quad \frac{u+v+1}u\\
\frac{u+v}v,&\quad -\frac uv
\endaligned
$$
Together with those
obtained by the action of $\tau$ these pairs $f\wedge (1-f)$ span a
lattice of rank $8$ in $\wedge^2(\Q(X)^*)$, which, luckily, contains
$x\wedge y$. These functions belong to the
multiplicative subgroup of $\Q(X)^*$ generated by $-1,u,v,u+1,v+1,u+v$
(modulo constants this subgroup consists precisely of  those functions
in $\Q(X)^*$ with divisor supported on $H$).
After a straightforward calculation we find that if we
let
$$
\aligned
z_1=& \frac{v^2}{u^2(u+1)}\\
z_2=&-\frac v{(u+1)^2}\\
z_3=&-v\\
z_4=&\frac{v(v+1)}{u(u+1)^2}\\
z_5=&-(u+v)\\
z_6=&-\frac{(v+1)}u\\
z_7=&\frac{u+v}v\\
\endaligned
$$
then
$$
3(x\wedge y)=3\langle z_1\rangle +\langle z_2\rangle
+\langle z_3\rangle+\langle z_4\rangle+\langle z_5\rangle
+2\langle z_6\rangle+2\langle z_7\rangle
$$
and therefore
$$
3\eta(x,y)=D(\xi),
$$
where
$$
\xi(w)=3[z_1]+[z_2]+[z_3]+[z_4]+[z_5]+2[z_6]+2[z_7], \quad w\in X(\C).
$$
There are of course many different choices of $z_j's$ that would give
a similar decomposition; as we will see below, our choice was
determined to simplify the calculation of $\m(A)$.

 We now would like to compute $m(A)$ for $A$ given
in (26) using Theorem 1. The calculation of $\partial \gamma$
is somewhat more complicated than in out previous examples since $A$
is reciprocal 
$$
A^*(x,y)=x^3y^2A(1/x,1/y)=A(x,y)
$$
so we cannot proceed as in \S 5.  However, by considering
the vanishing of $\disc_x A$, we
find that $\gamma$ can be chosen so that
$$
\partial \gamma
=[w_1]-[\overline w_1], \qquad w_1=\left(-2,\frac{1-\sqrt{-15}}2\right) \in
X(\Qbar).
$$
The two points $w_1$ and $\overline w_1$ map to the double point
$(x,y)=(1,-1)$ on $Y$. The hypotheses of Theorem~3 are met and
hence we find that
$$
\m(A) = r d_{15}
\tag 27
$$
for some $r\in \Q^*$. Numerically $r=1/6$;  we now show how we can
actually prove that $r = 1/6$ in this case.

Let $K=\Q(\sqrt{-15})$ and $\Gamma=PSL_2(\o_K)$. The group $\Gamma$
acts discretely on hyperbolic space $\h^3$; let $M=\Gamma \backslash \h^3$
be the corresponding orbifold quotient. By the theorem of Humbert
quoted in the introduction we have
$$
\Vol(M)=\pi d_{15}/6,
$$
where $\chi$ is the Dirichlet character attached to $K/\Q$. On the
other hand, Gangl [Ga] has shown that six times a fundamental
domain for the action of $\Gamma$ can be triangulated with
into ideal tetrahedra with shape parameters $2\delta$, where
$$
\delta=3[a]+4[(a+1)/2]+2[a+1]+2[a-1]\in \B(K), \qquad
a=(1+\sqrt{-15})/2,
\tag 28
$$
 which implies that
$$
\Vol(M)=\tfrac13 D(\delta).
$$
Computing the value of $\xi(w)$ at the boundary point $w_1$ we find that
$$
\xi(w_1)=3[(a+3)/4]+2[a-1]+2[a+1]+4[a/2].
\tag 29
$$
It is now easy to check that $\xi(w_1)$ and $\delta$ are the same
modulo the threefold symmetry $z,1/(1-z),1-1/z$ and hence represent
the same element of $\B(K)$. Putting everything together we find that
indeed $r=1/6$.

An independent proof that $m(A) = d_{15}/6$ is given by Dunfield in
the Appendix to this paper.  He constructs an explicit arithmetic
manifold whose A-polynomial is exactly our polynomial $A(x,y)$  The
formula then follows from Humbert's theorem, without the need for
constructing an explicit triangulation.

\endexample

\example{Example 4 --  Involutions and reciprocal models of an elliptic curve}
We now consider some examples where we start with the
elliptic curve and search for functions $x,y,z_1,\cdots,z_n$ on it
such that {\bf (T)} holds and such that $x,y$ are related by a reciprocal
polynomial relation $A(x,y)=0$.     The point here is that the mapping
$(x,y) \to (1/x,1/y)$ corresponds to an involution on the curve and that
involutions on elliptic curves are well understood.   In particular whether
or not the involution has a fixed point has an important effect on the 
formula for the Mahler measure of $A(x,y)$ as we illustrate below.

Consider
$$
X: \qquad v^2+uv+v=u^3-u
\tag 30
$$
 a minimal Weierstrass model of an elliptic curve of conductor
$14$. The rational points of $X$ consist of the group $G$ of order $6$
generated by $Q=(1,0)$. Consider $U\subset \Q(X)^*/\Q^*$ the subgroup
of functions, modulo constants, with divisors supported on $G$.

It is not hard to find a basis for $U$; for example,
$u,v,1-u,1+u,u+v$. The corresponding divisors are
$$
\matrix
\format \r\quad & \c\qquad & \r\quad &\r\quad & \r\quad &\r  \quad&\r\\
n & nQ & u & v &  1-u & 1+u & u+v\\
\\
0 & O & -2 & -3 & -2 & -2 & -3 \\
1 & (1,0) & 0 & 1 & 1 & 0 & 0\\
2 & (0,0) & 1 & 1 & 0 & 0& 3\\
3& (-1,0) & 0 & 1 & 0 & 2 & 0\\
4& (0,-1) & 1 & 0 & 0 & 0 & 0\\
5& (1,-2) & 0 & 0 & 1 & 0 & 0 \\
\endmatrix
\tag 31
$$

To ensure that  $A$ is reciprocal we pick an involution $\sigma$ of
$X$ defined over $\Q$ preserving $G$ and then we chose $x,y \in
\Q^*(X)$ among the functions $f$ satisfying
$$
f^\sigma=f^{-1}
\tag 32
$$
There are three separate cases, which we now describe.

a) Let
$$
\sigma: \quad P \mapsto -P+4Q.
$$
On the function field $\Q(X)$ it acts by
$$
u\mapsto -\frac{u+v}{u^2}\qquad v\mapsto -\frac{v(u+v)}{u^3}
$$
and on $G$ it gives  the permutation $\sigma=(04)(13)$, where the
number $n$ corresponds to the point $nQ \in G$. The divisor
$$
\delta=\sum_{n=0}^5c_n[nQ]
$$
 is the divisor of a function satisfying (53) if and only if
$$
c_{\sigma n}=-c_n, \qquad \sum_{n=0}^5c_n=0, \qquad
\sum_{n=0}^5nc_n=0;
\tag 33
$$
that is,
$$
\delta=a([0]-[4Q])+b([Q]-[3Q]), \qquad a+2b\equiv 0 \bmod 3.
$$
 The lattice of solutions to (33) is generated by $(a,b)=-(3,0)$ and
$(a,b)=-(1,1)$ corresponding to the functions
$$
x=v+1, \qquad y = \frac{u(u+1)}v
$$
respectively.

To compute $A$ we  eliminate $u$ and $v$ from the equations
$$
\aligned
0&=v^2+uv+v-u^3+u\\
0&=yv-u(u+1)\\
0&=x-v-1
\endaligned
$$
by taking successive resultants, factoring the result and
selecting the correct factor. We obtain
$$
A=(1-x)y^3+2y^2+2y-x+x^2
$$
with matrix of coefficients
$$
\matrix
\format \r\quad & \r\quad & \r\\
&-1&1\\
&2&\\
&2&\\
1&-1&\\
\endmatrix
$$
Note that $A$ is tempered, something we knew in advance as $\{x,y\}$ is
trivial. In fact, $A$ turns out to be the $A$-polynomial of the
hyperbolic $3$-manifold m009.

It is not too hard to compute a triangulation of $x\wedge y$ and we
find that we may take
$$
\xi = \left[-u\right] + \left[\frac{u+v}{u^2}\right] +
\left[\frac{v+1}{u^2}\right].
$$
The fixed points of the involution are those points $P$ on $E$
satisfying $2P=4Q$. The only two fixed points mapping to the torus
$|x|=|y|=1$ are $P_0$ and $\bar P_0$ where
$$
P_0=\left(-\tfrac12(1+\sqrt{-7}),-2\right),
$$
 which map to $x=-1,y=1$.  We have
$$
\xi(P_0)=2\left[\tfrac12(1+\sqrt{-7})\right]+
\left[\tfrac18(3+\sqrt{-7})\right],
$$
which, as it happens, are exactly the shapes that Snap gives for a
triangulation of $m009$. We conclude that
$$
\pi \m(A)=\Vol(m009)=\D(\xi(P_0)).
$$
As in Example 3, we can prove that $m(A) =d_7/2$ using Humbert's
formula for the volume, since $m009$ is an arithmetic manifold.

The fact that the conductor $14$ of $X$ and the discriminant $7$ of
the field $F = \Bbb Q(\sqrt{-7})$
both involve $7$ is not a coincidence.  Note that $F$ is the field of
definition of $P_0$ which is a point of finite order on $X$.

b) Repeating the above calculations for the involution
$$
\sigma: \qquad P \mapsto -P+Q
$$
we find we may take
$$
x=\frac{v^3}{u(u-1)^2(u+1)},\qquad y=-\frac{(u+1)(u+v)}{u^2(u-1)}
$$
related by  $A(x,y)=0$ with $A$ given by
$$
\matrix
\format \r\quad & \r\quad  & \r\quad & \r\quad & \r\\
&&&1&\\
&1&1&1&1\\
&5&-4&5&\\
1&1&1&1&\\
&1&&&\\
\endmatrix
$$
 This new polynomial in turn is the $A$-polynomial of the hyperbolic
$3$-manifold m221. The only fixed point of $\sigma$ on the torus is
$x=1,y=-1$ which corresponds to the four points $P$ on $E$ satisfying
$2P=Q$. These four points are Galois conjugates of each other and are
defined over a quartic field with $r_1=1$, $\disc=-2^6\cdot7$ and
defining equation $x^4-2x^3+x^2-2x+1$. Geometrically the point
$x=1,y=-1$ corresponds to the complete structure on $m221$ and we again
find that
$$
\pi \m(A) = \Vol(m221) = 24 V_F,
$$
where, as usual, the value $24$ here denotes a rational number equal
to $24$ to 40 decimal places (and we know that $\pi\m(A)\sim V_F$)

Note that examples a) and b) show the same curve arising as the zero
locus of the $A$-polynomial of two non-isomorphic hyperbolic
manifolds

c) Finally, we consider the fixed point free involution
$$
\sigma: \qquad P \mapsto P + 3Q.
$$
We find that
$$
x=\frac v{u(u+1)}, \qquad y = \frac v{1-u}
$$
are inverted by $\sigma$ and are related by the equation $A(x,-y)=0$
where $A$ is
$$
\matrix
\format \r\quad & \r\quad & \r\\
1&1&\\
1&1&1\\
&1&1\\
\endmatrix
$$
In this case the path $\gamma$ is closed; its projection via $x,y$
does not cross the torus, though it does have a part lying entirely on
it. The symbol $\{x,y\}$ is {\it not} trivial; in fact,
$\m(A)=.2274812230$ is not zero and, as it turns out, is the
smallest measure of a two variable integer polynomial known
\cite {BM}. According to the conjectures of Bloch--Beilinson $m(A)$
should be commensurable with $L'(E,0)$ and, numerically, it does seem
to be the case that $\m(A)=L'(E,0)$ \cite {D, BEM, RV1}.

\endexample

\example{Example 5 -- The A-polynomial of m389}
In this example, we consider the A-polynomial of the
one cusped hyperbolic manifold
$m389$, the complement of the knot $k5_{22}$ of \cite{CDW}, another name
for the 10-crossing
knot $10_{139}$ of \cite{R}.   The interest here is that we can show that
$\pi m(A) = \Vol{m389}$ and hence obtain a formula for $m(A)$ in terms of the
Dedekind zeta function of the invariant trace field of m389.     Another interesting
feature is that the curve defined by $A(x,y) = 0$ is of genus 1 and the method
of computing $A(x,y)$ leads directly to an explicit triangulation of $x \wedge y$.

The manifold m389 can be triangulated with 5 ideal
tetrahedra with shape parameters $\xi_1, \dots, \xi_5$   The complex number
$\xi_j$ is the cross-ratio of the points in $\Bbb C$ representing the vertices
of the $j$th tetrahedron, using
the upper half space model of hyperbolic 3-space $\Bbb H^3$  \cite {Th}.

To compute the A-polynomial we start from the
deformed gluing equations for the manifold:
$$z_1 z_2 z_3 z_4 (1-z_2)^{-1}(1-z_5) = -1 $$
$$z_2^{-1}z_4^{-1} z_5^{-1}(1-z_1)^{-1}(1-z_2)(1-z_3)^{-1}(1-z_4)(1-z_5) = -1$$
$$z_3 z_5 (1-z_1)(1-z_2)^{-1}(1-z_4) = 1$$
$$z_1^{-1} z_2 (1-z_1)(1-z_3)^{-1}(1-z_4)^{-1}(1-z_5)^{-1}=-1 $$
$$z_2^{-1}z_3^{-2}(1-z_1)^{-1}(1-z_2)(1-z_3)^2 (1-z_4)^{-1}(1-z_5)^{-1}= -1$$
$$z_1^{-1}z_3^{-1}z_4^{-1}z_5(1-z_1)(1-z_2)(1-z_4)^{-1}(1-z_5)^{-2}=M^2$$
$$z_3 z_4^2 z_5^{-1}(1-z_1)(1-z_2)^{-1}(1-z_4)^{-1}(1-z_5)^2=L^2$$
The first 5 of these equations express the condition that the dihedral angles of the
tetrahedra meeting at an edge of the triangulation add up to $2\pi$.  These
are exactly as obtained from the program Snap.   The last two equations
given by Snap can be written as $\mu = 1$, and $\lambda = 1$ which state
that the holonomy around a meridian and longitude of the cusp are each
integer multiples of $2\pi$.    We replace $\lambda$ by $\lambda_1=\lambda/\mu$
for a reason given below and then deform the equations to $\mu=M^2$ and
$\lambda_1=L^2$ where $L$ and $M$ are ``deformation parameters'', obtaining
the sixth and seventh equations given above.

Snap also gives us an exact solution of these equations
$z_j = \xi_j$, $L=M=1$, with the
$\xi_j$ lying in the shape field $F$ of the manifold, a field
 of type $4, [2,1], -688$ with defining polynomial $x^4 - 2x - 1$.

The A-polynomial is obtained by successively eliminating the variables
$z_1,\dots,z_5$ from the 7 equations, thus obtaining a
polynomial $A(L,M^2)$.   This polynomial is even in $M$ since the manifold
is a knot complement and we have chosen the longitude to have trivial mod 2 homology
\cite {CCGLS}.
At each stage of the process of elimination, we retain only equations that vanish at
$(z_1,\dots,z_5,L,M) = (\xi_1,\dots,\xi_5,1,1)$. The polynomial $A(L,M^2)$ so
obtained is an irreducible (over $\bar \Bbb Q$)
factor of the A-polynomial defined in \cite{CCGLS}.  In general $A(x,y)$ will
have coefficients in a finite extension of $\Bbb Q$ but usually they are
rational integers, as in this example.
By a slight abuse of terminology, we will call $A(x,y)$ {\it the} A-polynomial
of $m389$.

For our example, the coefficient matrix of $A$ is
$$
\matrix
\format \r\quad & \r\quad & \r\quad & \r\quad &\r\\
  & 1 &   &   &   \\
  &-3 & 0 &-1 &   \\
1 & 7 & 6 & 7 & 1 \\
  &-1 & 0 &-3 &   \\
  &   &   & 1 &
\endmatrix
$$

Using Maple's {\tt algcurve} package, one finds that $A(x,y) = 0$ is of
genus $1$ and that defining
$$x = -(u+1)v/(u+v+1),$$
$$y = u^5(u+v+1)v^{-2}(v+1)^{-2},$$
reduces this to the Weierstrass form
$$X:  v(v+1) = u^2(u+1),$$
which is the curve $43A1$ in Cremona's tables \cite{Cr}.   The curve
is of rank $1$ with no torsion and $X(\Bbb Q)$ is generated by the
point $Q =[0,0]$ (in $[u,v]$ coordinates).   The involution
$(x,y) \to (1/x,1/y)$  is given by $P \to Q-P$.

The main fact that we need from geometry is that
the triangulation of the symbol
$x \wedge y$ is given by $\sum_{j=1}^5 z_j \wedge (1-z_j)$.
This is ultimately a consequence of Schl\"afli's formula for change in
volume of a hyperbolic polyhedron under deformation of its
edges and dihedral angles \cite {Ho}, \cite {Du2}, \cite {CCGLS}, 
\cite{NZ}.

Using the system of equations obtained in the construction of
$A(x,y)$, and the relation between $(x,y)$ and $(u,v)$, we can
recursively solve for the $z_j$ and obtain
$z_1 = z_5 = (u+v+1)/(u+1)$, $z_2 = 1/(u+1)$, $z_3 = (v+1)/(u+v+1)$,
$z_4 = (u+1)/u$ giving the explicit triangulation
$$x \wedge y = \langle z_1 \rangle + \langle z_2 \rangle +
\langle z_3 \rangle + \langle z_4 \rangle + \langle z_5 \rangle.$$

We find that
$\disc_y(A(x,y)) = y^3 f(y) (y-1)^{12}$, where $f(y)$ does not
vanish on $|y|=1$ and that $A(x,1) = (x+1)^4$.   The point
$(x,y)=(-1,1)$ is the only singularity of the curve
$Y = \{A(x,y)=0\}$ and the points $(u,v)$ of $X$ lying above $(-1,1)$
are the four Galois conjugates of $(\alpha,\beta)$ , where
$\alpha$ is a root of the quartic $x^4 - 2x - 1$ defining $F$
and $\beta = \a^3 - 1$
The corresponding values of $z_1,\dots,z_5$ are the Galois
conjugates of $\a^3-1,(-\a^3+\a^2-\a+3)/2,(\a^3-1)/2,\a^3-1,\a^3-1$.
One of these is the geometric solution $\xi$ of the gluing equations with
$\IM \xi_j > 0$ for all $j$ and this is found to correspond to the
conjugate of $\alpha$, say $\a^{(1)}$ with $\IM \a^{(1)} < 0$.
We let $\a^{(2)}$ be its complex conjugate and $\a^{(3)}$,
$\a^{(4)}$ denote the real conjugates.   The corresponding
values of  $D(\xi^{(k)})$ are $\Vol(m389),-\Vol(m389),0,0$,
for $k=1,\dots,4$.

Near $y = 1$,
the four roots of $A(x,y) = 0$ have Puiseux expansions
$x_k(y) = -1 + c^{(k)}(y+1) + O((y+1)^2)$, with $c^{(k)} \in F^{(k)}$,
where the $F^{(k)}$ are the four embeddings of $F$ into $\Bbb C$.
Choose the numbering so that $c^{(1)}$ and $c^{(2)}$ are complex
conjugates with $\IM c^{(1)} > 0$, and $c^{(3)}$, $c^{(4)}$
are real.   Then it is not hard to see that $|x_1(y)| > 1 > |x_2(y)|$
for all $|y| = 1$ with $y \ne 1$, and that $|x_3(y)| = |x_4(y)| = 1$
for all $|y| = 1$  \cite{BoH}.
Hence  $2\pi m(A) = \int_\gamma \eta(x,y)$,
where the endpoints of $\gamma$ are the two non-real points
$(\a^{(1)},\b^{(1)})$ and $(\a^{(2)},\b^{(2)})$ of E lying
above $(-1,1)$.    We thus see that
$$2\pi m(A) = \Vol(m389)-(-\Vol(m389)) = 2\Vol(m389),$$
and thus
$$\pi m(A) = \Vol(m389).$$
(It is clear that the signs
must be as indicated since both $m(A)$ and $\Vol(m389)$ are positive.)

By Borel's theorem, we know that $\Vol(m389) = r V_F$ for some
rational $r$.
Numerically, $\Vol(m389) = 4.8511707573$, from Snap, and
$V_F = 0.4042642297$ from Pari, so $r = 12$
to 28 decimal places.

\endexample

\example {Example 6 -- The A-polynomial of v2824 -- a harmless $\times$-crossing} In most of our
examples, the contributions to $m(A)$ have come from one or two points
of the form $(\pm 1, \pm 1)$.  Here we describe the A-polynomial of
the manifold $v2824$ where there is a contribution from another point
on the torus.  This is a one-cusped manifold triangulated by 7
tetrahedra with shapes in a field of type $8, [6,1], -397538359$.
Another field that comes into the picture is the shape field of
$v2824(1,0)$ which Snap finds to be of type $9, [5,2],
4497953501$.  By a computation as in the previous example, we
find that the A-polynomial is of the form $A(L,M^2)$ where $A(x,y)$ is
of degree $20$ in $x$ and $11$ in $y$, and has height $228$ (so we do
not display it here).  We find that $\disc_x(A) = -F_1 F_2^2 F_3^6
(y-1)^{146} y^{78}$, where $F_1(y)$ is of degree $40$ and does not
vanish on $|y|=1$, $F_2(y)$ has 5 pairs of complex zeros on $|y|=1$
but none corresponds to a crossing of $|x| = 1$ by a root $x_k(y)$ of
$A(x,y) = 0$.  The polynomial $F_3(y) = y^2+y+1$ vanishes at the cube
roots of unity $\z$ and $\bar \z$ and the corresponding solutions of
$A(x,y)=0$ are $(\z^2,\z)$ and its complex conjugate.

This point corresponds to what we call
an {\it $\times$-crossing} of the circle $|x| = 1$, i.e. such that the corresponding
solutions $x_k(y)$ with $x_k(\z) = \z^2$ cross the circle with
a non-zero and non-infinite slope.  Thus $(\z^2,\z)$ is the endpoint
of two separate paths in the cycle $\gamma$ of (6).

One verifies that this is an $\times$-crossing by examining
the Puiseux expansion at $(\z^2,\z)$ and finds that
$$A(\z^2 + st , \z+s) = 27q(t)s^3 + O(s^4),$$
where $q(t) = 167t^3 - 21(1+\z)t^2 + 81\z t + 11$ defines
a cubic extension of $\Bbb Q(\z_3)$ of relative discriminant
$-23$.  In fact the splitting field of $q(t)$ over $\Bbb Q$ is
the compositum  $F = Q(\a,\z_3)$ where $\a$ solves $x^3-x-1$,
defining a field of type $3, [1,1], -23$.   So $F$ itself is
of type $6, [0,3], -3^3 \cdot 23^2$.

As in the previous example, the method of computation of $A(x,y)$
gives us the corresponding shapes for the solutions with
$(x,y) = (\z^2,\z)$.   These shapes $\xi_j$ lie in the sextic field $F$
and the six values of $D(\xi^\sigma)$ for the corresponding six Galois
conjugates of $\xi$ are $0,0,c,-c,-c,c$,  where
$c = .94270736$.   The contribution to $\pi m(A)$ from the point
$(\z^2,\z)$ is thus $2c = 1.88541472$.    By Theorem~2,
since $(x,y)_w$ is clearly trivial at $w=(\z^2,\z)$ this $\xi \in B(F)$.
Formula (29) does not immediately apply here since $r_2(F) = 3$.
However, if $G$ is the cubic subfield $\Bbb Q(\a)$ of $F$, we compute
$V_G = D(\alpha)/12 = .078558946$ and numerically find that
$c = D(\xi) = 12 V_G$ to $40$ decimal places.  The fact
that $D(\xi) \sim V_G$ is a consequence of Theorem B.

The final formula for $\pm m(A)$ involves 5 terms: the two just described,
apparently giving $12 V_G$,  a contribution from $(-1,1)$, giving
$\Vol(v2824)=6.0463301699$, which according to our Theorem is a rational
multiple of $V_H$ for $H$ the shape field of $v2824$ described above, and
numerically we have $\Vol(v2824) = (71/6)V_H$; the final two terms are the
absolute values of the components of the Borel regulator of $v2824(1,0)$.
and are given by Snap to be
$17443417807$ and $1.8283473415$. These come from the
two branches above $(1,1)$.  For the definition of the Borel
regulator, see \cite {NY, CGHN}.

\endexample

\example{Example 7 -- The A-polynomial of m410 -- a less benign $\times$-crossing}

In this example, we illustrate that the condition {\bf (B)} of
Theorem~2 need not always hold, even for A-polynomials
of hyperbolic manifolds.
The manifold $m410$ is triangulated by 5 equal tetrahedra of
shape $\z_3 = (1 + \sqrt{-3})/2$ and hence has volume
$5 D(\z_3) = 5 \pi d_3$.
If we interchange the roles of meridian and longitude as given
by SnapPea, then
this manifold  has an A-polynomial $A(L,M^2)$
where $A(x,y)$ has the Newton polygon

$$
\matrix
\format \r\quad &\r\quad &\r\quad &\r\quad & \r\quad
          &\r\quad &\r\quad & \r\quad & \r\quad & \r\quad &\r\\
   & 1  & -2 &  3 & -1 & -1 &  3 & -2 &  1 &    &    \\
-1 &  4 & -1 & -2 &  3 &-10 &  3 & -2 & -1 &  4 & -1 \\
   &    &  1 & -2 &  3 & -1 & -1 &  3 & -2 &  1 &
\endmatrix
$$

We check that $\disc_x(A) = F_1 F_2^2 (y-1)^2 y^5$, where $F_1(y)$
does not vanish on $|y| = 1$ and $F_2(y) = y^6 + y^5 + y^3 + y + 1$
has 2 pairs of roots $\a,\bar \a, \b, \bar \b$ on $|y| = 1$ with
arguments $\pm .736837\cdots$ and $\pm 1.72698\cdots$.  There is an
$\times$-crossing at $\b, \bar \b$ but not at $\a, \bar \a$, so, as in the
previous example, the point $(x_0,y_0)=(1-\b+\b^2-\b^3+\b^5,\b)$ on
$A(x,y) = 0$ and its complex conjugate are endpoints of a path in the cycle
$\gamma$ of (6).  Solving for the corresponding shape vector $\xi$ as
above, we find that the shapes lie in field $F$ of type $12, [0,6],
2^8 5^2 17^2 37^4$, and that the 12 values of $D(\xi)$ are
$[0,0,0,0,-a,a,-a,a,b,-b,b,-b]$ where $a = 3.85787307$ and $b =
0.30130595$.  The value $a$ does not contribute to $\pi m(A)$
since it corresponds to a point $(x_1,y_1)$ on $A(x,y) = 0$ which does
not lie on the torus.  Thus $\pi m(A) = 5 D(\z_3) + 2 b$.  In this
case, unlike the previous example, $\xi$ does not lie in the Bloch
group $B(F)$ since it is easy to see that $x_0$ and $y_0$ are
multiplicatively independent.

\endexample

\example{Example 8 -- Chinburg's conjecture for $d_{19}$, $d_{40}$ and $d_{120}$}

The polynomial
$A(x,y) = (x+x^2+x^3+x^4+x^5)(y^2+1)+(1+6x+2x^2-8x^3+2x^4+6x^5+x^6)y$ with
Newton polygon
$$
\matrix
 \format \r\quad &\r\quad &\r\quad & \r\quad & \r\quad & \r\quad &\r\\
    &  1  &  1  &  1  &  1  &  1  &  \\
  1 &  6  &  2  & -8  &  2  &  6  &  1  \\
    &  1  &  1  &  1  &  1  &  1  &
\endmatrix
$$
was found by examining a large number of polynomials of this general
shape.   The curve $A(x,y)=0$ has a single singularity at $(1,1)$ and
the corresponding field of the normalization is $\Bbb Q(\sqrt{-19})$,
as revealed by the Puiseux expansion at $(1,1)$.    The curve is
in fact the elliptic curve 475C1 of conductor $5^2 \cdot 19$.
Numerical calculation of $m(A)$ shows that $m(A) = (2/5)d_{19}$ holds
to 50 decimal places.

We can prove that $m(A) = r d_{19}$ for some $r \in \Q^*$. To see this
note that the involution $\sigma: x\mapsto 1/x, y \mapsto 1/y$ on $E$
has fixed points and hence $E/\langle \sigma \rangle$ is isomorphic to
$\p^1$ and by Galois descent $\{x,y\}$ is torsion. Since $\partial
\gamma$ maps to $(1,1)$ via $(x,y)$ our claim follows from the main
theorem.   This adds another to the list of discriminants for which
we have an answer to the conjecture of Chinburg \cite {Ch, BRV}.

Two further discriminants we can add are $-40$ and $-120$.  For the 
discriminant $40$, if $A(x,y) = (x^2-x+1)(x^2+1)(y^2+1) + 
x(14x^2-32x+14)y$ then numerically $m(A) = d_{40}/6$ and if
$A(x,y) = (x^2+1)(x^2+x+1)(y+1)^2 - 24x^2y$, then numerically
$m(A) = d_{120}/36$.  In each case $A(x,y) = 0$ is
of genus 0 with a single singularity at $(1,1)$ which gives the only
contribution to $m(A)$, and thus Theorem~3 implies
that $m(A) \sim d_f$, for $f = -40, -120$, respectively.

These are just a few of many different examples of the general
form $A(x,y) = B(x)(y^2+1) + C(x)y$ for which
one can prove that $m(A) = rd_f + sd_g$ for rational numbers $r,s$.
Typically both $r$ and $s$ are non-zero so these do not necessarily
contribute to Chinburg's conjecture, but occasionally there is only
a single term.   For example, the polynomial with
$B(x) = (x^2+x+1)(x^4-x^2+1)$ and $C(x) = -x^6+14x^5-32x^3+14x-1$
has $m(A) = rd_3 + sd_{219}$  where presumably $r = 1/2$ and $s = 1/72$.

Another interesting example has $B(x) = (x^4+x^3+x^2+x+1)^2$ and
$C(x) = 2x^8+4x^7-x^6-17x^5-26x^4-17x^3-x^2+4x+2$.   Remarkably
$A(x,y) = 0$ is of genus $1$ and in fact is the elliptic curve
15A3 of conductor $15$.  We find that
$$\disc_y(A) = - (7x^2+11x+7)(x^2+x+1)(x-1)^2(x+1)^2(2x^2+3x+2)^2 x^2,$$
so there are singularities on the torus with
$$x = -1, 1, (-1\pm\sqrt{-3})/2, (-3\pm\sqrt{-7})/4, (-11\pm 5\sqrt{-3})/14.$$
In spite of this, only the singularities at $(-1,-1)$ and $(1,1)$ 
contribute to $m(A)$ giving  $m(A) = 4d_3/5 + 26d_3/5 = 6d_3$, a rather 
dull result for such an interesting polynomial.    Here the rational 
number $6$ has been
surmised from numerical computations but in this case one would expect
to be able to prove this rigorously without recourse to geometry.

\endexample

\example{Example 9 -- the A-polynomial of $v1859$ and other manifolds}

The manifold $v1859$ is an arithmetic manifold with invariant trace
field $Q(i)$ and volume $6 D(i)$, with $i = \sqrt{-1}$.    It has the
very simple A-polynomial $1 + iL - iM - LM$ in the standard $(L,M)$
coordinates.   Since $A$ is even in neither $L$ nor $M$ the relation
one would expect is $2\pi m(A) = \Vol(v1859) = 6D(i)$.   However, it
is easy to see that $2\pi m(A) = 4D(i)$  \cite{Sm, BEM},
so $2\pi m(A) = (2/3)\Vol(v1859)$.

Nathan Dunfield explains that this is in
fact just one of a sequence of manifolds obtainable by Dehn surgery
on one of the cusps of a manifold $\Cal N \Cal R$ constructed by Neumann
and Reid  \cite {NR}.   From SnapPea, one obtains $\Cal N \Cal R$
 as the unique
2-cover of $m135$ with $H_1 = \Z/2 + \Z/2 + \Z + \Z$.   This
is a two cusped manifold with {\it strong geometric isolation},
which means that Dehn surgery on one cusp does not affect the shape
of the other cusp.   This means that the A-polynomials of all the
manifolds $\Cal N \Cal R(a,b,0,0)$ are the same,
namely $1 + iL - iM - LM$,
the same as that of the manifold $m135$ \cite {Du1}.

Now consider the sequence of one cusped manifolds
$\M_n = \Cal N \Cal R(n,1,0,0)$.
By the results of Thurston \cite {Th}, $\Vol(\M_n)$ is a strictly increasing
sequence bounded below by $\Vol(m135) = 4D(i)$ and above by
$\Vol(\Cal N \Cal R) = 2\Vol(m135) = 8D(i)$.  But $2\pi m(A) = 4D(i)$ is
constant and hence never equal to $\Vol(\M_n)$.   Our example,
$v1859$ turns out to be $\M_2$ and the fact that it is arithmetic
is somewhat accidental.   The manifolds $\M_n$ are not arithmetic
for $n > 2$ and their volumes are not commensurable with $D(i)$.

\endexample

\subheading{7. The Universal Triangulation}
 It is possible to construct, in some sense, the
{\it universal triangulation} of a symbol $\{x,y\}$ as follows. Let
$u_1,u_2, \ldots, u_{n+1}$ be independent variables. Given $m_1, m_2,
\ldots, m_{2n} \in \Z^{n+1}$ consider the subvariety of $(\C^*)^{n+1}$
defined by the system of equations
$$
X: \qquad \left \{
\aligned
u^{m_1}+u^{m_2}&=1\\
u^{m_3}+u^{m_4}&=1\\
\vdots \\
u^{m_{2n-1}}+u^{m_{2n}}&=1
\endaligned
\right .
$$
where we use the multinomial notation
$u^l=u_1^{l_1}\cdots u_{n+1}^{l_{n+1}}$ for
$l=(l_1,l_2,\ldots,l_{n+1})\in \Z^{n+1}$.
It will be convenient to switch to matrix notation. Let
$M=(m_1,m_2,\ldots, m_{2n}) \in \Z^{(n+1)\times 2n}$ be the 
integral matrix with columns $m_1,m_2,\ldots, m_{2n}$. Then the
condition
$$
u_1\wedge u_2 = \sum_{k=1}^n u^{m_{2k-1}}\wedge u^{m_{2k}}
$$
 in the free abelian group generated by the $u$'s is equivalent to
the matrix equation
$$
MJM^T=
\pmatrix
0&1\\ -1&0
\endpmatrix
\perp 0 \perp \cdots \perp 0
\tag{34}
$$
where
$$
J= \pmatrix
0&1\\ -1&0
\endpmatrix
\perp \cdots \perp \pmatrix
0&1\\ -1&0
\endpmatrix
$$
is the standard symplectic form \cite {NZ}. 
If this holds then on the variety $X$
we have the triangulation
$$
u_1\wedge u_2 = \langle u^{m_1}\rangle+\cdots + \langle u^{m_{2n-1}}\rangle.
$$
 Since we have $n+1$ variables and $n$ equations the variety $X$
will generically be of dimension $1$; hence, to every generic matrix
$M$ satisfying (34) we can associate a curve $X$ and a
pair of rational functions $u_1,u_2$ on it with trivial symbol
$\{u_1,u_2\}\in K_2(X)$.

Note that since we only want the triangulation up to torsion we can
just as well consider equations defining $X$ of the form $\pm u^m\pm
u^{m'}=1$ instead of $u^m +u^{m'}=1$.

 We may want to construct in this way polynomials that have all the
known characteristics of the $A$-polynomial of a hyperbolic
$3$-manifold. For this we would need to ensure in addition that the
projection of $X$ onto the $u_1,u_2$ coordinates be invariant under
the involution $\sigma: u_1\mapsto 1/u_1,u_2\mapsto 1/u_2$ so that the
resulting polynomial cutting out this image is reciprocal. This is
somewhat more complicated and we will content ourselves with
considering the following example.

Let
$$
M=
\pmatrix
-1&0&1&0&0&0\\
0&0&0&1&1&0\\
1&0&1&0&0&1\\
0&1&0&1&2&0
\endpmatrix
$$
 It is a simple matter to check that $M$ satisfies (34). For
convenience, we will label the variables as $x,y,u,v$ instead of
$u_1,\ldots,u_4$. Consider the variety $X$ determined by this matrix
with the following choice of signs
$$
\left\{
\aligned
ux^{-1}+v&=1\\
ux+vy&=1\\
-v^2y-u&=1
\endaligned
\right.
$$
Using resultants or Macaulay2 it is not hard to eliminate $u$ and $v$
from this system and find that the projection of $X$ onto the $x,y$
coordinates is cut out by the polynomial $A(x,y)$
$$
\matrix
\format \r\quad & \r\quad & \r\quad & \r\\
 & & &1\\
1&-2&-2&1\\
1& & & \\
\endmatrix
$$
 which is the $A$-polynomial of the manifold m009 (and it was in fact
starting backwards from this polynomial and a triangulation for this
manifold that we were led to consider the matrix $M$).

The involution $\sigma: x\mapsto 1/x,y\mapsto 1/y$ on the projection
lifts to $(\C^*)^4$ preserving $X$ via $u\mapsto u, v\mapsto vy$.

We now consider $M'=MU$ for $U\in Sp_6(\Z)$. Clearly, $M'$ satisfies
(34) since by definition $UJU^T=J$; we would also like to have
$\sigma$, as an involution of the ambient torus $(\C^*)^4$, preserve
$X$. It is easy to check that this will hold if $U$ has the form
$$
\pmatrix
\pmatrix
a&b\\c&d
\endpmatrix
&&\\
& \pmatrix
a&b\\c&d
\endpmatrix
&\\
&&
\pmatrix
1&0\\0&1
\endpmatrix
\endpmatrix
$$
for any $\gamma = \pmatrix
a&b\\c&d
\endpmatrix \in SL_2(\Z)$.
We obtain in this way infinitely many curves $X$ together with a
trivial symbol $\{x,y\}$ in $K_2(X)$ whose minimal polynomial $A(x,y)$
is reciprocal. The corresponding volume function $V$ on $X$ with
$dV=\eta(x,y)$ is
$$
V=D((u/x)^av^c)+D((ux)^a(vy)^c)+D(-v^2y).
$$
We give some examples below.

The involution $\sigma$ always has at least two fixed points which are
zeroes of $A$, namely $(\pm 1, 1)$; the corresponding values of $v$,
from which we get the value of $u=-1-v^2$, satisfy the equations
$$
(\pm 1,1): \qquad (\mp (1+v^2))^av^c+(\mp (1+v^2))^bv^d=1
$$
respectively. Typically there are no other points that contribute to
the value of $m(A)$.

\bigskip

a) For  $\gamma=  \pmatrix
-1&0\\0&-1
\endpmatrix
$
we have the system of equations
$$
\left\{
\aligned
u^{-1}x+v^{-1}&=1\\
u^{-1}x^{-1}+v^{-1}y^{-1}&=1\\
-v^2y-u&=1
\endaligned
\right .
$$
and the corresponding polynomial $A$ is
$$
\matrix
\format \r\quad & \r\quad &\r\quad &\r\quad & \r\quad & \r\quad & \r\\
 & & & & &1& \\
 &1&1&-4&-2&0&1\\
-1&0&2&4&-1&-1& \\
 &-1& & & & & \\
\endmatrix
$$
 As it turns out, this polynomial is the $A$ polynomial of the
hyperbolic $3$-manifold m367. The points $(\pm
1, -1)$ are the only intersections of $A(x,y)=0$ with the real torus
$|x|=|y|=1$ and the corresponding values of $v$ satisfy the cubic
equations $v^3-v^2+2v-1=0$ and $v^3-v^2-1=0$  respectively. Let
$F_1$ and $F_2$ be the fields these equations determine; they have
discriminants $-23$ and $-31$ respectively and both have
$r_2=1$. Therefore, the hypotheses of Theorem~3 apply and we
conclude that $\pi \m(A)$ is a rational linear combination of
$V_{F_1}$ and $V_{F_2}$.

\bigskip

b) For $\gamma= \pmatrix
-1&1\\0&-1
\endpmatrix
$
we find $A$ to be
$$
\matrix
\format \r\quad & \r\quad &\r\quad &\r\quad &  \r\quad & \r\quad &
\r\quad & \r\quad & \r\\
 & & & &1& & & & \\
 & & &1&-3&-2&-2&1&1\\
 &-1&-2&1&10&1&-2&-1& \\
1&1&-2&-2&-3&1& & & \\
 & & & &1& & & & \\
\endmatrix
$$
Again the only intersection points of $A(x,y)=0$ and the real torus
are $(\pm 1,1)$ and the corresponding points in $X$ lie in quartic
fields with $r_2=2$ and discriminants $117$ and $229$ respectively. We
do not know if $A$ is the $A$-polynomial of a hyperbolic manifold in
this case.

\bigskip

c) For  $\gamma= \pmatrix
1&1\\0&1
\endpmatrix
$
we find $A$ to be
$$
\matrix
\format \r\quad & \r\quad &\r\quad &\r\quad & \r\quad & \r\\
-1& & & & & \\
 &0&3&1&-1& \\
 &1&-1&-3&0& \\
 & & & & &1\\
\endmatrix
$$
which turns out to be the $A$-polynomial of the hyperbolic manifold
s254.

Again the only relevant points are $(\pm 1,1)$ resolving in $X$ over
a cubic field $v^3+v^2+v+2=0$ of discriminant $-83$ and the quadratic
field $\Q(\sqrt{-3})$. Theorem~3 applies and hence $\pi \m(A)$
is a rational linear combinations of the invariant $V_F$ of these two
fields.

\subheading{8. Final remarks}

1) One may ask if all reciprocal polynomials in $\A$ are $A$-polynomials of
   $1$-cusped manifolds. It is clear that they parameterize the
   deformation of a collection of hyperbolic tetrahedra for which $V$
   is the resulting volume function but what is not entirely clear is
   how these tetrahedra are glued together and what the resulting
   geometric object might be.

2) One could ask if for the $A$-polynomial $A$ of a hyperbolic
   manifold $M$ there is a topological interpretation, possibly in
   terms of representations of $\pi_1(M)$, for the points of $\partial
   \gamma$ contributing to the value of $\pi\m(A)$. (Note, for
   example, that condition {\bf B} of Theorem~3 precisely
   corresponds to Dehn surgery points in the representation variety.)
   Such an interpretation would further link the geometry of $M$ and
   the Mahler measure of its $A$-polynomial.

3) On the negative side, however, we should point out that Example~9
   shows that $\m(A)$ and $\Vol(M)$ can be completely unrelated.

4) We do not believe that there are any continuous families such as
   those considered in \cite{BEM, RV1} of
   polynomials in $\A$ (the tame symbol at finite primes should give a
   non-trivial obstruction for the symbol $\{x,y\}$ to be trivial).

5) Note that all the examples of polynomials $A \in \A(\Q)$ that we
   have exhibited have integer coefficients (with $A$ normalized to
   have vertex coefficients $\pm 1$). This is likely not to be a
   coincidence; in general, we expect that if $\{x,y\}$ is in
   $K_2({\Cal X})$ for ${\Cal X}$ a regular model of the curve $X$
   then the minimal polynomial $A$ of $\{x,y\}$ should have integer
   coefficients and conversely. This fact should be analogous to the
   case of $K_1({\Cal O}_F)={\Cal O}_F^*$ with ${\Cal O}_F$ the ring
   of integers of a number field $F$, which is trivial to verify.


\head
{\bf Appendix: The Mahler measure of the $A$-polynomial of $\bold m129(0,3)$}
\endhead

\head Nathan M. Dunfield \endhead

This appendix gives an alternate proof of the exact value of the
Mahler measure of the polynomial $A$ of \S6, Example 3.  It turns out
that there is a hyperbolic 3-manifold $N$ whose $A$-polynomial, an
invariant of the $\SL{2}{\C}$-representation variety of $\pi_1(N)$, is
essentially $A$.  For the $A$-polynomial of a manifold, the Mahler
measure is often related to the volume of the manifold itself.  In the
case of $N$, one has $\pi m(A) = \Vol(N)$.  The manifold $N$ shares a
finite cover with the Bianchi orbifold $B = \PSLn{2}{(\o_{-15})}\backslash\h^3$
and so its volume is rationally related to that of $B$; in fact the
two volumes are equal.  As discussed above, a formula of Humbert's
says that $\Vol(B) = \pi d_{15}/6$.  Putting this together shows that
$m(A) = d_{15}/6$.  This approach is more than just a replacement for
the geometric input from Gangl's theorem used in
Example~3---it removes the need to find the $z_i$'s at all.

The contents of this appendix are summarized as follows.
Section~\Aa\  gives the definition of the $A$-polynomial.
Section~\Ab\  gives an alternate formulation of the $A$-polynomial, and
clarifies the relationship between the standard $A$-polynomial and the
one coming from the gluing equation variety which is used in
the body of this paper.  Section~\Ab\ does not relate directly to the
computation of $m(A)$, and can be skipped by uninterested readers.
Section~\Ac\ describes in detail the computation of the Mahler measure of
$A$ following the outline above.

\subheading{Acknowledgments} This work was done while the author was
at Harvard University, and was partially supported by an NSF
postdoctoral fellowship.

\heading{\Aa. The $A$-polynomial of a 3-manifold}\endheading

Consider a compact orientable 3-manifold $M$ whose boundary, $\partial M$,
is a torus.  We are interested in the typical case where the interior
$M - \partial M$ has a complete hyperbolic metric of finite volume; that is
$M - \partial M = \Gamma \backslash\h^3$ where $\Gamma$ a non-uniform lattice in $\Isom^{+}(\h^3)
= \PSLn{2}{\C}$.  In this case, a neighborhood of the deleted $\partial M$ is
a finite-volume torus cusp.  Such an $M$ is referred to as a 1-cusped
hyperbolic 3-manifold.

The $A$-polynomial of $M$ is an invariant of the space of
representations of $\pi_1(M)$ into $\SL{2}{\C}$, as viewed from $\partial M$.
It was introduced in \cite{CCGLS}.  The target group $\SL{2}{\C}$ is used
because it is almost the group of isometries of $\h^3$; the
identification of $\pi_1(M)$ with a lattice in $\PSLn{2}{\C}$ gives an
interesting representation $\pi_1(M) \to \PSLn{2}{\C}$, called the
holonomy representation.  The reason that $\SL{2}{\C}$ is used instead
of $\PSLn{2}{\C}$ is mostly just historical; initially, $\SL{2}{\C}$ was
chosen because it can be easier to work with.

To begin, consider a finitely generated group $\Gamma$ and let $G$ be
either $\SL{2}{\C}$ or $\PSLn{2}{\C}$.  Let $R(\Gamma, G)$ be the set
of all representations $\Gamma \to G$.  The set $R(\Gamma, G)$ has a
canonical structure as a complex algebraic variety.  This can be
defined by choosing generators $\gamma_1, \ldots , \gamma_n$ of
$\Gamma$ and embedding $R(\Gamma, G) \to G^n$ via $\rho \mapsto
(\rho(\gamma_1), \ldots, \rho (\gamma_n))$.  It is natural to study
representations up to inner automorphisms of $G$, that is, to consider
$X(\Gamma, G) = R(\Gamma, G)/G$ where $G$ acts via conjugation.
Technically, one has to take the algebro-geometric quotient to deal
with orbits which are not closed; this way, $X(\Gamma, G)$ is also a
complex algebraic variety.  The original space $R(\Gamma,G)$ is called
the representation variety of $\Gamma$, and the quotient $X(\Gamma,
G)$ is called the character variety.  If $Y$ is a topological space,
set $R(Y,G) = R(\pi_1(Y), G)$ and $X(Y, G) = X(\pi_1(Y),G)$.

For a 1-cusped hyperbolic 3-manifold, the character variety $X(M, G)$
has deep connections to both the geometry and topology of $M$  First,
there is the irreducible component $X_0(M, G)$ of $X(M, G)$ which
contains the equivalence class $[\rho_0]$ of the holonomy representation
coming from the hyperbolic structure.  Thurston showed that the
(complex) dimension of $X_0(M, G)$ is always one, giving a curve of
distinct representations.  Moreover, $X(M,G)$ is in fact an affine
variety, so $X_0(M,G)$ is non-compact.  It turns out that if one adds
on ``ideal points'' to $X_0(M,G)$ to compactify it, these ideal points
correspond to certain topologically essential surfaces in $M$ itself.
For more background and applications of character varieties to the
study of 3-manifolds, see the survey \cite{Sh}.

To define the $A$-polynomial, we first need to understand the
character variety of the torus $\partial M$.  The fundamental group of
$\partial M$ is just $\Z \times \Z$.  Geometrically, a pair of generators
$(\mu, \lambda)$ of $\pi_1(M)$ corresponds to a pair of simple closed curves
in $\partial M$ which meet in a single point.  The generators $(\mu, \lambda)$
are usually called the meridian and longitude respectively.  Since
$\pi_1(M)$ is commutative, any representation $\rho \maps \pi_1(\partial M) \to
\SL{2}{\C}$ is reducible, and the corresponding M\"obius transformations
have a common fixed point on $P^1(\C)$.  Moreover, if no element of
$\rho(\pi_1(\partial M))$ is parabolic, $\rho$ is conjugate to a diagonal
representation with
$$
\rho(\mu) = \left( \matrix M & 0 \\ 0 & M^{-1} \endmatrix
\right) \quad \mbox{and} \quad \rho(\lambda) = \left( \matrix L & 0 \\ 0 &
    L^{-1} \endmatrix \right).
$$
As such, $X(\partial M, \SL{2}{\C})$ is approximately $\C^* \times \C^*$ with
coordinates being the eigenvalues $(M,L)$.  This isn't quite right, as
switching $(M, L)$ with $(M^{-1}, L^{-1})$ gives a conjugate
representation.  In fact, $X( \partial M, \SL{2}{\C})$ is exactly the
quotient of $\C^* \times \C^*$ under the involution $(M, L) \mapsto (M^{-1},
L^{-1})$

Now the inclusion $i \maps \partial M \to M$ induces a regular map $i^* \maps
X(M, \SL{2}{\C}) \to X(\partial M, \SL{2}{\C})$ via restriction of
representations from $\pi_1(M)$ to $\pi_1(\partial M)$.  Let $V$ be the
1-dimensional part of the $i^*\big(X(M, \SL{2}{\C})\big)$.  More precisely,
take $V$ to be the union of the 1-dimensional $i^*(X)$, where $X$ is
an irreducible component of $X(M, \SL{2}{\C})$.  The curve $V$ is used
to define the $A$-polynomial.  To simplify things, we look at the
plane curve $\overline{V}(M, \SL{2}{\C})$ which is inverse image of
$V$ under the quotient map $\C^* \times \C^* \to X(\partial M, \SL{2}{\C})$.  The
$A$-polynomial is the defining equation for $\overline{V}(M,
\SL{2}{\C})$; it is a polynomial in the variables $M, L$.  Since all
the maps involved are defined over $\Q$, the $A$-polynomial can be
normalized to have integral coefficients.

\heading{\Ab . Additional comments on the $A$-polynomial}\endheading

In this section, I discuss an alternate definition of the
$A$-polynomial, the definition of the $\PSLn{2}{\C}$ variant of the
$A$-polynomial, and relationship between the character variety and the
gluing equation variety.  Readers interested primarily in the
computation of the Mahler measure of the polynomial in Example 3
should skip ahead to Section~\Ac .

\subheading{\Ab .1.~Alternate Definition}  
There is an equivalent formulation of the $A$-pol\-y\-no\-mi\-al which is
useful in understanding the $\PSLn{2}{\C}$ version of the
$A$-polynomial and its relationship to the gluing equation variety
(see Section~2.3).  Let $\overline{R}(M, \SL{2}(\C))$ be the subvariety
of $R(M, \SL{2}{\C}) \times P^1(\C)$ consisting of pairs $(\rho, z)$ where
$z$ is a fixed point of $\rho(\pi_1(\partial M))$.  Let $\overline{X}(M,
\SL{2}{\C})$ be the algebro-geometric quotient of $\overline{R}(M,
\SL{2}{\C})$ under the diagonal action of $\SL{2}{\C}$ by conjugation
and M\"obius transformations respectively.  There is a natural regular
map $\pi \maps \overline{X}(M, \SL{2}{\C}) \to X(M, \SL{2}{\C})$ which
forgets the second factor.  The inverse image $[ \rho ]$ under $\pi$ is
isomorphic to the subset of $P^1(\C)$ fixed by $\rho(\pi_1(\partial M))$; this
is $2$ points unless either $\rho(\pi_1(\partial M))$ contains a parabolic (in
which case it is $1$ point) or is trivial (in which case it is
$P^1(\C)$).  In particular, for the geometric component $X_0(M,
\SL{2}{\C})$, the map $\pi$ has fibers generically consisting of $2$
points.

The advantage of the augmented character variety $\overline{X}(M,
\SL{2}{\C})$ is that given $\gamma \in \pi_1(\partial M)$ there is a regular
function $e_\gamma$ which sends $[ (\rho, z) ]$ to the eigenvalue of
$\rho(\gamma)$ corresponding to $z$.  In contrast, on $X(M, \SL{2}{\C})$
only the trace of $[\rho(\gamma)]$ is well-defined.  Moreover, the map
$\overline{X}(\partial M, \SL{2}{\C}) \to \C^* \times \C^*$ given by $e_\mu \times e_\lambda$
is an isomorphism.  The $1$-dimensional part of the image of $i^*
\maps \overline{X}(M, \SL{2}{\C}) \to \overline{X}(\partial M, \SL{2}{\C}) =
\C^* \times \C^*$ is again the plane curve $\overline{V}(M, \SL{2}{\C})$
whose defining equation is the $A$-polynomial.

\subheading{\Ab .2.~The $A$-polynomial over $\PSLn{2}{\C}$} As I'll
describe, the alternate definition of the $A$-pol\-y\-no\-mi\-al given
above easily adapts to the $\PSLn{2}{\C}$ version of the
$A$-polynomial described in \cite{Cha}.  One difference which
introduces a slight technicality, though, is that you
gain one representation of $\Z \times \Z$ by passing to
$\PSLn{2}{\C}$.  In particular, there is a representation $\sigma$
whose image is $\Z/2 \times \Z/2$, and which acts on $P^1(\C)$ {\it
  without\/} a common fixed point.  This corresponds to picking two
geodesics in $\h^3$ that meet in a right angle and taking the group
generated by rotation by $\pi$ about each of them.  The representation
$\sigma$ is the only representation of $\Z \times \Z$ that does not
lift to $\SL{2}{\C}$, and it cannot be deformed to any other
representation.  Thus $X(\partial M, \PSLn{2}{\C})$ is reducible, with
a 2-dimensional component which is essentially the same as in the
$\SL{2}{\C}$ case, and an isolated point consisting of $[\sigma]$.

To define the $\PSLn{2}{\C}$ $A$-polynomial, first consider the subset
$R'(M, \PSLn{2}{\C})$ of $R(M, \PSLn{2}{\C})$ consisting of
representation which do {\it not\/} restrict to $\sigma$ on
$\pi_1(\partial M)$.  Then,
just as before, we can define the augmented character variety
$\overline{X}(M, \PSLn{2}{\C})$ as the quotient of the appropriate
subset of $R'(M, \PSLn{2}{\C}) \times P^1(\C)$.  If $\gamma \in
\pi_1(\partial M)$, the eigenvalue function $e_\gamma$ is no longer
well-defined; however, its square the \emph{holonomy function}
$h_\gamma$, still is.  In particular, the value of function $h_\gamma
\maps \overline{X}(M, \PSLn{2}{\C}) \to \C^*$ at $[ (\rho, z) ]$ is
the derivative of $\rho(\gamma)$ at the fixed point $z$.  The pair of
functions $h_\mu \times h_\lambda$ gives an isomorphism of
$\overline{X}(\partial M, \PSLn{2}{\C})$ with $\C^* \times \C^*$.
Analogously, denote by $\overline{V}(M, \PSLn{2}{\C})$ the
$1$-dimensional part of the image of $i^* \maps \overline{X}(M,
\PSLn{2}{\C}) \to \overline{X}(\partial M, \PSLn{2}{\C})$.  The
defining equation of this plane curve is the $\PSLn{2}{\C}$
$A$-polynomial.
 
It's natural to ask what the differences are between these two
versions of the $A$-polynomial.  Consider the map $p \maps
\overline{V}(M, \SL{2}{\C}) \to \overline{V}(M, \PSLn{2}{\C})$ coming
from restricting the projection $\overline{X}(\partial M, \SL{2}{\C}) \to
\overline{X}(\partial M, \PSLn{2}{\C})$.  The important point is that $p$ may
not be onto $\overline{V}(M, \PSLn{2}{\C})$.  This is because not every
representation $\rho \maps \pi_1(M) \to \PSLn{2}{\C}$ lifts to
$\SL{2}{\C}$.  However, the holonomy representation coming from the
hyperbolic structure always lifts \cite{Cu}, and so there is a
always a relationship between the ``geometric components'' of the two
plane curves.  Also, the obstruction to $\rho$ lifting lies in $H^2(M,
\Z/2)$, so if this group vanishes (as it does for the complement of a
knot in $S^3$) every representation lifts.
 
\subheading{\Ab .3.~Gluing equations and the character variety}
There is another way of looking at the character variety, namely by
considering solutions to the ``gluing equations'' of a triangulation,
as is done in \S6, Example 5.  For completeness, I'll outline the
precise relationship between this point of view and the standard one;
for a more detailed account see \cite{Cha}.
Throughout, the reader unfamiliar with the gluing equations can find
additional details in \cite{Th, Ch. 4}  or
alternatively, in \cite{NZ, BP, Ra}.

Fix an ideal triangulation of the 1-cusped manifold $M$.  That is, a
choose a cell-complex $\mathscr{T}$ where all the cells are
tetrahedra, and where the complement of the single $0$-cell is
homeomorphic to $M - \partial M$.  Given such a triangulation, we can try
to build a hyperbolic structure on $M$ by assigning hyperbolic shapes
to each of the (purely topological) tetrahedra in $\mathscr{T}$.  The
shape of an oriented ideal geodesic tetrahedron $\Delta$ in $\h^3$ is
described by a complex number $z$.  More precisely, this \emph{shape
  parameter} $z$ is associated to an edge of $\Delta$.  Any
tetrahedron isomorphic to one is with vertices at $\{ 0, 1, \infty, z \}$
under the action of $\PSLn{2}{\C}$; for such a tetrahedron the
invariant of the edge with vertices $\{ 0, 1 \}$ is $z$.  The shape
parameter is also the cross ratio of 4 vertices in a suitable order.
Since the four vertices are distinct, $z$ is in $P^1(\C) - \{0,1,\infty\}$.

 For any tetrahedron $\Delta$, opposite
edges of $\Delta$ have the same parameter.  The three distinct edge parameters
$z_1, z_2, z_3$ of $\Delta$ are related by $z_2 = 1/(1-z_1)$ and $z_3 =(z_1
- 1)/z_1$.  Another way of saying this is that the possible triples
of shapes are exactly the solutions in $\C^3$ to the equations
$$
z_1 z_2 z_3 = -1 \quad \mbox{and} \quad  z_1 z_2 - z_1 + 1 = 0.
\tag{\ShapeEqn}
$$
Notice that the solutions in $\C^3$ to these equations do not include
any of the degenerate shapes $\{0, 1, \infty\}$.

Returning to our triangulation $\mathscr{T}$, suppose it consists of
tetrahedra $\Delta_1, \ldots, \Delta_n$.  An assignment of shapes to these
tetrahedra is given by a point in $(\C^3)^n$ which satisfies $n$
copies of the shape equations (\ShapeEqn ).  The \emph{gluing
  equation variety}, $G(\mathscr{T})$, is the subvariety of possible
shapes where we require in addition that the following \emph{edge
  equations} are satisfied: for each edge in $\mathscr{T}$, the
product of the edge parameters of the tetrahedra around that edge is
$1$.  This requirement says that the hyperbolic structures on the
individual tetrahedra glue up along the edge.

So what exactly does $G(\mathscr{T})$ parameterize?  Consider the
universal cover $\widetilde{M}$ of $M - \partial M$ with induced ideal
triangulation $\widetilde{\mathscr{T}}$.  As I'll explain, points in
$G(\mathscr{T})$ correspond to functions $\widetilde{M} \to \h^3$ which
take every ideal tetrahedra in $\widetilde{\mathscr{T}}$ to a
non-degenerate geodesic ideal tetrahedron, and which satisfy the
following: given a deck transformation $g \in \pi_1(M)$, there exists an
isometry $\tau$ of $\h^3$ such that $f \circ g = \tau \circ f$.  For such an
$f$, there is an associated holonomy map $\rho$ defined by $g \mapsto \tau$;
then $f$ is a pseudo-developing map for $\rho$ in the sense of
\cite{Du2, \S2.5}.  How does a point $p \in G(\mathscr{T})$
correspond to such a map?  Well, you start with some $\Delta_0$ in
$\widetilde{\mathscr{T}}$ and send it to a geodesic ideal tetrahedra
in $\h^3$ with the specified shape.  Then you extend the map to an
adjacent tetrahedra and repeat.  The edge equations exactly ensure
that this process of ``analytic continuation'' is well defined.  In
fact, it's not hard to see that the points in $G(\mathscr{T})$ are in
bijective correspondence with such maps $\widetilde{M} \to \h^3$,
provided that we regard two maps as the same if they differ by an
orientation preserving isometry of the target.

Because of the holonomy representation associated a point in
$G(\mathscr{T})$ we have a function
$$
G(\mathscr{T}) \to X(M, \PSLn{2}{\C}),
$$
which can be shown to be regular.  Moreover, there is a regular map
from $G(\mathscr{T})$ to the augmented character variety
$\overline{X}(M, \PSLn{2}{\C})$ for the following reason.  Consider a
collar neighborhood $N$ of $\partial M$ which, in the cell complex
$\mathscr{T}$, is just a small neighborhood of the unique $0$-cell.
In each tetrahedra, $N$ is a small neighborhood of the four vertices
cut off by four triangles.  Let $\widetilde{N}$ be the inverse image
of $N - \partial M$ in $\tilde{M}$.  The subgroup $\pi_1(\partial M)$
consists of the deck transformations which stabilize some particular
connected component $\widetilde{N}_0$ of $\widetilde{N}$.  Consider
the pseudo-developing map $f \maps \widetilde{M} \to \h^3$ coming from
a point in $G(\mathscr{T})$, and let $\rho$ be its holonomy
representation.  Then $f(\widetilde{N}_0)$ is made up of standard
pieces of ideal tetrahedra, and its closure intersects $\partial \h^3
= P^1(\C)$ in a single point $z$.  The point $z$ must be a fixed point
of $\rho (\pi_1(\partial M))$, and this is exactly the additional
information need to construct the (regular) map $G(\mathscr{T}) \to
\overline{X}(M, \PSLn{2}{\C})$.

One can show that the map $G(\mathscr{T}) \to \overline{X}(M,
\PSLn{2}{\C})$ is injective.  However, it need not be surjective.  The
problem is that while every representation $\pi_1(M) \to
\PSLn{2}{\C}$ has some kind of reasonable pseudo-developing map
$\widetilde{M} \to \h^3$, you can't always make it geodesic with
respect to a fixed triangulation of $M$.  In particular, when you to
try straighten things out, some edges may shrink off to $\partial\h^3$.
While you can see \cite{Du2, \S2.5} for more
details, let me give a simple example.  

Let $M$ be the exterior of the figure-8 knot in $S^3$.  It has a
standard ideal triangulation $\mathscr{T}_0$ with two tetrahedra.  For
this triangulation, $G(\mathscr{T}_0)$ is an irreducible curve
\cite{Tu, Ch.~4}.  However, $ \overline{X}(M,
\PSLn{2}{\C})$ has two components; the geometric one and one consisting
of reducible representations.  It's the latter component that doesn't
appear in $G(\mathscr{T}_0)$.  In addition, $M$ has another
triangulation $\mathscr{T}_1$ with 5 tetrahedra which has an edge of
valence one.  In this case, $G(\mathscr{T}_1) = \emptyset$ because one of the
edge equations is $z_i = 1$ and the shape equations (\ShapeEqn )
force $z_i$ not to be $1$.

Finally, for any triangulation, it is easy to write down holonomy
functions $h_\gamma \maps G(\mathscr{T}) \to \C^*$ in terms of the shapes of
the tetrahedra.  For a basis $(\mu, \lambda)$ of $\pi_1(\partial M)$, the map $h_\mu
\times h_\lambda \maps G(\mathscr{T}) \to \C^* \times \C^*$ gives another plane curve,
which is a union of components of $\overline{V}(M, \PSLn{2}{\C})$; its
defining equation is thus a factor of the $\PSLn{2}{\C}$
$A$-polynomial (see \cite{Cha} for more).

\heading{\Ac. Computation of $m(A)$ using the $A$-polynomial} \endheading

\subheading{\Ac .1.~Definition of $N$} 
I'll begin with an explicit description of $N$.  While I said above
that $N$ is a manifold, in fact it is an orbifold with one cusp and a
singular locus consisting of a knot labeled $\Z/3\Z$; however, the
$A$-polynomial makes just as much sense in this context.  The orbifold
$N$ will be described as Dehn filling on \SnapPea\ census-manifold
$m129$ \cite{W1, HW}.  The hyperbolic manifold
$m129$ has two cusps and is the complement of the Whitehead link in
$S^3$  The manifold $N$ is the $(0, 3)$ orbifold Dehn filling on the
first cusp of $m129$.  Here, I'm using the standard \SnapPea\ basis for
the homology of the cusps of $m129$, which is the basis which makes
the cusp shape square; this differs from the one you would choose if
you regard $m129$ as the complement of the Whitehead link.  With
respect to the latter choice of basis, $N$ is the $(6, 3)$ Dehn
filling, with the orientation convention that the $(1,1)$ Dehn surgery
on the Whitehead link gives the trefoil, not the figure-8.

The orbifold $N$ is hyperbolic and its volume is about
$3.13861389446$.  This can be confirmed with \Snap\ \cite{CGHN}, which
verifies the existence of a hyperbolic structure using exact
arithmetic in $\Q(\sqrt{-3}, \sqrt{5})$.  This also follows from the
argument below which shows that $N$ shares a common cover with the
Bianchi orbifold $B$; the implication is not immediate, however, as it
uses Thurston's Orbifold Theorem to say that $N$ must have {\it some\/}
geometric decomposition.

\subheading{\Ac .2.~How $N$ was found} The polynomial $A$ of Example~3
satisfies all the known consequences of being the $A$-polynomial of a
manifold, and so it was natural to search for a manifold with that
polynomial.  As you'll see in \S \Ac .5, when you try to compute the
Mahler measure of $A$, you get contributions from only 2 points.
This suggests that if $A$ is the $A$-polynomial of a manifold, then
the volume of that manifold should be a rational multiple of $\pi
m(A)$.

The orbifold $N$ was found by a brute-force search of orbifolds which
appeared numerically to have volume which was a rational multiple of
$\pi m(A)$.  First, I looked in the complete \SnapPea\ census of
cusped hyperbolic manifolds with fewer than 8 tetrahedra
\cite{HW, CHW}, but didn't find
anything.  I then took those census manifolds which had more than one
cusp, and tried various Dehn fillings to get some additional
candidates.  Motivated by the fact that the Bianchi group for
$\o_{-15}$ has 2 and 3 torsion, I allowed Dehn fillings which gave
orbifolds as well as manifolds.  This resulted in a small number of
possible examples.  In addition to $N$, I also looked at $M =
m412(3,0)$.  While $M$ doesn't have $A$-polynomial equal to $A$, is
also related to the Bianchi orbifold $B$ and will be used in the proof
that $N$ and $B$ share a common finite cover.

\subheading{\Ac .3.~Relationship to $\PSLn{2}{(\o_{-15})}$} Next, I'll
show that $N$ and the Bianchi orbifold $B = \PSLn{2}{(\o_{-15})}
\backslash \h^3$ are \emph{commensurable}, that is, have a common
finite cover.  I will verify this using the program \SnapPea\ 
\cite{W1}.  This is rigorous because \SnapPea\ never reports that two
manifolds are isometric without finding triangulations of the two
manifolds which are combinatorially isomorphic---in particular, this
involves only integer computations.  The converse, showing manifolds
are \emph{not} isometric is not rigorous in \SnapPea , because it uses
floating point computations to compute the canonical triangulations
\cite{W2, SW} which it then compares; for large triangulations,
round-off error can result in the wrong triangulation being identified
of as the ``canonical'' one.

Let $B$ be the Bianchi orbifold for $\o_{-15}$.  Baker showed that $B$
is Dehn surgery on a link in $S^3$ \cite{Bak}, and this description
was used to describe this orbifold to \SnapPea .  Let $M$ be the
orbifold $m412(3,0)$.  First, I note that $M$ and $B$ share a common
6-fold cover which is characterized by the fact that it is a manifold
with 2-cusps, has homology $\Z/3 + \Z/9 + \Z^2$, and is an irregular cover
of $M$ and a cyclic cover of $B$.  Thus $M$ and $B$ are commensurable
and $\Vol(M) = \Vol(B)$.  Now I claim $M$ and $N$ are commensurable
and that $\Vol(M) = \Vol(N)$.  In fact, they have a common 3-fold
cover which is a manifold with 2-cusps and homology $\Z/6 + \Z^2$.
This manifold is an irregular cover of both $M$ and $N$.  Now
commensurability is an equivalence relation, so therefore $N$ and $B$
are commensurable, and, moreover, $\Vol(N) = \Vol(B)$.  Combining this
with Humbert's formula for the volume of $B$, we have $\Vol(N) = \pi
d_{15}/6$.

\subheading{\Ac .4~Computing the $A$-polynomial of $N$}
In this section, I will outline the computation that the $A$-polynomial
of $N$ is
$$
A_N = 
(M^2 L +1) 
(M^6 L^2-M^4 L^3-3 M^4 L^2-M^4 L+M^2 L^2+3 M^2 L+M^2-L).
\tag{\ApolyOfN}
$$
I'll denote the second factor, which is the geometric one, by
$A'_N$.  This polynomial is closely related to $A$, and in fact $A'_N
= M^{-2} A( -M^2 L, -M^2)$; in particular $m(A'_N) = m(A)$.

To compute the $A$-polynomial of $N$, I will use the following
description of the fundamental group of $m129$:
$$
\pi_1(m129) = \spandef{a,b,c}{A c^2 A C b C = [a, b] = 1}
$$
where here $A = a^{-1}$, etc.  With this presentation, a
meridian-longitude basis for the fundamental group of the first cusp
of $m129$ is $\{ cA, acBAc \}$ and a basis for the second cusp is
$\{a,b\}$.  (The simplicity of the basis for the second cusp is why we
choose this presentation, rather than a standard 2-generator
presentation.)  Now $N$ is obtained by doing $(0,3)$ Dehn filling on
the first cusp of $m129$.  So $\pi_1(N)$ is gotten from
$\pi_1(m129)$ by adding the requirement that the longitude $\lambda =
acBAc$ has order 3, i.e.
$$
\pi_1(N) = \spandef{a,b,c}{A c^2 A C b C = [a,b] = \lambda^3 = 1}.
$$

To compute the $A$-polynomial of $N$, we take the
$\SL{2}{\C}$-representation variety of $\pi_1(N)$ and project it onto
the representation variety of $\pi_1(\partial N) = \left\langle a,b \right\rangle$.
Since $a$ and $b$ commute, their images under a representation $\rho
\maps \pi_1(N) \to \SL{2}{\C}$ have a common fixed point, and so $\rho$ is
conjugate into the form:
$$
\rho(a) = \left(\matrix a_0 & a_1 \\ 0 & a_3 \endmatrix \right) , \quad
\rho(b) = \left(\matrix b_0 & b_1 \\ 0 & b_3 \endmatrix \right) , 
\quad \mbox{and} \quad
\rho(c) = \left(\matrix c_0 & c_1 \\ c_2 & c_3 \endmatrix \right).
$$
Let $V$ be the variety of all representations of $\pi_1(N)$ into
$\SL{2}{\C}$ which have this form.  The variety $V$ splits into two
components, depending on whether $\rho(\lambda)$ has order 1 or 3.
The interesting case is the component $V_0$ where $\rho(\lambda)$ has
order 3.  The defining equations of $V$ are
$$
 \det \rho(a) = 1, \det \rho(b) = 1,
\det \rho(c) = 1, \rho(A c^2 A) = \rho(c B c), \tr \rho(\lambda) = -1,
$$
where the last equation is equivalent to $\lambda$ having order 3.
To get the contribution of $V_0$ to the $A$-polynomial, we project
$V_0$ onto the eigenvalues of $\{ \rho(a), \rho(b) \}$, that is, onto
the coordinates $\{a_0, b_0\}$.  Since we will be interested in the
1-dimensional components of the image, we can just project those
representations where $\rho(\left< a, b \right>)$ contains no
parabolics; this allows us to restrict to the subvariety $V'_0$ where
$a_1 = b_1 = 0$.  Using Gr\"obner bases in Macaulay 2 \cite{GS}, it is
easy to compute that the defining polynomial of the projection of
$V'_0$ is
$$
(M^2-L) (M^6 L^2-M^4 L^3-3 M^4 L^2-M^4 L+M^2 L^2+3 M^2 L+M^2-L).
$$
One can show that the first factor $M^2 - L$ comes from the subvariety of
reducible representations; by convention, it is not included in the
$A$-polynomial of $N$.  If one looks at the other half of $V$ where
$\rho(\lambda) = 1$, one gets a contribution to the $A$-polynomial of
$M^2 L+1$.  This second component is just the $A$-polynomial of
$m129(0,1)$, which is Seifert fibered.  Putting this together, we have
the formula (\ApolyOfN) for the $A$-polynomial of $N$.

\subheading{\Ac.5.~The Mahler measure of $A'_N$}
In this final section, I'll show that $\pi m(A'_N) = \Vol(N)$.
Combining this with $\Vol(N) = \pi d_{15}/6$ will give that $m(A'_N) =
d_{15}/6$, completing the calculation of the Mahler measure of the
polynomial in \S6, Example~3.  The close connection here between the
Mahler measure of the $A$-polynomial and the volume of the manifold is
common, though not universal; this type of connection was discovered
by Boyd \cite{BH}.  Roughly, the reason is that the 1-form $\eta$
of \S2 used to compute the Mahler measure has another meaning on the
curve $\overline{V}(N, \SL{2}{\C})$: it measures the change in volume
of representations $\rho \maps \pi_1(N) \to \SL{2}{\C}$.  In our very
simple case, the integral expression (4) for the Mahler measure boils
down to the difference between the volume of the holonomy
representation  and its complex conjugate.   

I'll now outline the computation for our particular $A'_N$; for
further details see \cite{BH}.  It is a little easier to work with
$$
B(x,y) = -y^2 x^3 + (y^3 - 3 y^2 + y) x^2 + (-y^2 + 3 y - 1) x + y,
$$
where $A'_N = B(L, M^2)$.  Write $B$ as $(-y^2) \prod_{k =
  1}^3\big(x - x_k(y)\big)$.  We are interested in the functions $x_k$
on the torus, where they can be made continuous.  For notational
convenience, set $x_k(t) = x_k(e^{i t})$ for $t \in [0, 2 \pi]$.
Because $B$ is reciprocal, one has that, after relabeling, $|x_1(t)| =
1/|x_2(t)|$ and $|x_3(t)| = 1$.  As all of the $|x_k(t)| = 1$ only for
$t = 0$ or $2 \pi$, we can define our choice of branches by $|x_1(t)| >
1$, $|x_2(t)| < 1,$ and $|x_3(t)| = 1$ for $t \in (0,2 \pi)$.  Formula
(4) then gives
$$
m(B) = \sum_{k = 1}^3 \frac{1}{2 \pi} \int_{0}^{2 \pi} \log^+ |x_k(t)| \ dt = \frac{1}{2 \pi} \int_0^{2 \pi} \log |x_1(t)| \ dt,
$$
since $m(-y^2) = 0$ and  only $|x_1(t)|$ has non-zero $\log^+$.  

As described in \cite{BH, \S 4}, if $\rho_t \maps \pi_1(N) \to
\SL{2}{\C}$ for $t \in [0,2 \pi]$ is a family of representations where
$e_\lambda(\rho_t) = x_1(t)$ and $e_\mu(\rho_t) = e^{i t/2}$, then the
derivative of volume function of these representations is $- \log |
x_1(t) | \ dt$.  Thus $\pi m(B) = (1/2)\big(\Vol(\rho_0) - \Vol(\rho_{2
  \pi})\big)$. It is easy to check, following
\cite{BH, \S 4}, that $\rho_t$ can be chosen so that $\rho_0$ is the
holonomy representation of the hyperbolic structure on $N$, and $\rho_{2
  \pi}$ is the complex conjugate of $\rho_0$.  Thus as $\Vol(\rho_0) =
\Vol(N) = - \Vol(\rho_{2 \pi})$, we have $\pi m(A'_N) = \pi m(B) =
\Vol(N)$.

\enddocument

\end{thebibliography}

\end{document}